\documentclass[12pt,a4paper,final]{iopart}

\usepackage[toc,page]{appendix}

\usepackage{color}
\usepackage{iopams}  
\usepackage{graphicx}
\expandafter\let\csname equation*\endcsname\relax

\expandafter\let\csname endequation*\endcsname\relax

\usepackage{iopams} 
\usepackage{amssymb}
\usepackage{amsmath}
\usepackage{amsthm}
\usepackage{dsfont}
\usepackage{mathrsfs}
\usepackage{multirow}
\usepackage{tikz,pgfplots}
\usetikzlibrary{patterns}
\usepackage{subcaption}
\pgfplotsset{compat=newest}
\usepackage{float,indentfirst,algorithm}
\usepackage{pgfplots}
\usepackage{booktabs}
\usepackage{hyperref}
\usepackage{multicol}
\usepackage{nomencl}
\usepackage{algpseudocode}
\usepackage{ulem}
\makenomenclature
\makeatletter\let\list@fname\nomname \makeatother

\newtheorem{theorem}{Theorem}

\newtheorem{remark}{Remark}

\usepackage[capitalise]{cleveref}
\usepackage{appendix}
\usepackage{multicol}

\newcommand{\mean}{m}

\newcommand{\N}{\mathcal{N}}

\newcommand{\F}{\mathcal{F}}

\newcommand{\bigO}{\mathcal{O}}
\newcommand{\E}{\mathbb{E}}

\newcommand{\R}{\mathbb{R}}

\newcommand{\dd}{\mathrm{d}}

\newcommand{\bG}{\mathcal{F}}

\newcommand{\by}{x}

\newcommand{\Dt}{\Delta \tau}
\newcommand{\Range}{{\rm Range}}

\definecolor{darkred}{rgb}{.7,0,0}

\begin{document}

\title[Enhancing Exploration in Ensemble Kalman Inversion]{Enhancing Exploration in Ensemble Kalman Inversion via Covariance Inflation and Projected Dropout}

\author{Baojun Che}
\address{School of Mathematical Sciences, Peking University, Beijing, China}
\ead{bjche25@stu.pku.edu.cn}
  
\author{Daniel Zhengyu Huang}
\address{Corresponding author. Beijing International Center for Mathematical Research,  Center for Machine Learning Research, Peking University, Beijing, China}
\ead{huangdz@bicmr.pku.edu.cn}
  
\author{Xinying Mao}
\address{School of Mathematical Sciences, Peking University, Beijing, China}
\ead{maoxinying@stu.pku.edu.cn}
  
\author{Weijie Wang}
\address{School of Mathematical Sciences, Peking University, Beijing, China}
\ead{wwj66285509350@stu.pku.edu.cn}

\begin{abstract}
Ensemble Kalman inversion (EKI) is a derivative-free and parallelizable method for parameter estimation in Bayesian inverse problems involving expensive forward models. It adapts Kalman filtering by repeatedly assimilating the observation data through an artificial dynamical system with identity state evolution and the parameter-to-observation map as its observation operator.
However, its performance can be limited by ensemble collapse and confinement to the affine span of the initial ensemble.
We propose an exploration-enhanced EKI framework that combines multiplicative covariance inflation with a projected-dropout correction of the ensemble mean. Inflation helps maintain ensemble spread, while projected dropout enables the mean to explore directions beyond the initial affine span.
For linear inverse problems, we provide a discrete-time, finite-ensemble analysis of the resulting method and prove the least-squares objective converges exponentially to its minimum in expectation. Within the initial ensemble subspace, the contraction factor is independent of the conditioning of the inverse problem, and the limiting ensemble covariance provides a low-rank representation of the posterior covariance.
Experiments on ill-conditioned linear systems, 100-dimensional nonlinear benchmarks, and a PDE-constrained barotropic flow inverse problem demonstrate the effectiveness of the proposed method.
\end{abstract}

\vspace{2pc}

\submitto{\IP. The authors are in alphabetical order.}




\section{Introduction}

Many problems arising in scientific computing, data assimilation, and uncertainty quantification can be formulated as inverse problems, in which one seeks to recover an unknown parameter from indirect and noisy observations. A common formulation is the weighted nonlinear least-squares problem
\begin{equation}
\label{eq:intro-nlsq}
\min_{\theta\in\mathbb R^{N_\theta}}
\frac12
\bigl(x-\mathcal F(\theta)\bigr)^T
\Sigma_\nu^{-1}
\bigl(x-\mathcal F(\theta)\bigr)
=
\min_{\theta\in\mathbb R^{N_\theta}}
\frac12
\left\|
\Sigma_\nu^{-1/2}\bigl(x-\mathcal F(\theta)\bigr)
\right\|^2,
\end{equation}
where \(\theta\in\mathbb R^{N_\theta}\) is the unknown parameter, \(x\in\mathbb R^{N_x}\) is the observation, \(\mathcal F:\mathbb R^{N_\theta}\to\mathbb R^{N_x}\) is the forward map, and \(\Sigma_\nu\in\mathbb R^{N_x\times N_x}\) is symmetric positive definite. This framework covers both unregularized data-fitting problems and Tikhonov-regularized inverse problems.

In many large-scale applications, the parameter dimension is high, the forward model is computationally expensive, and derivative information with respect to \(\theta\) is either unavailable or unreliable. Ensemble Kalman inversion (EKI)~\cite{chen2012ensemble,emerick2013investigation,iglesias2013ensemble} provides a derivative-free and ensemble-based approach for approximately solving~\eqref{eq:intro-nlsq}. The method introduces the artificial filtering model
\begin{subequations}
\label{eq:intro_dynamics}
\begin{align}
  &\text{evolution:}    &&\theta_{n+1}=\theta_n, \label{eq:intro_evolve0}\\
  &\text{observation:}  &&x_{n+1}=\mathcal F(\theta_{n+1})+\nu_{n+1}, \label{eq:intro_obser0}
\end{align}
\end{subequations}
where \(\nu_{n+1}\sim\mathcal N(0,\Sigma_\nu)\), and applies ensemble Kalman filtering ideas~\cite{kalman1960new,evensen1994sequential} to the repeated observation \(x_{n+1}=x\). In this way, EKI inherits the computational simplicity, scalability, and parallel structure of ensemble Kalman methods while avoiding adjoint or gradient computations.

Despite these advantages, standard EKI has two structural limitations. First, the static artificial evolution in~\eqref{eq:intro_dynamics} does not replenish ensemble spread, while the Kalman analysis step continuously reduces covariance. This leads to ensemble collapse and algebraically slow convergence: in the linear setting, continuous-time EKI has \(L^2\) mean error of order \(\bigO(t^{-1/2})\)~\cite{schillings2017analysis,schillings2018convergence,blomker2019well}, and discrete-time EKI has corresponding order \(\bigO(n^{-1/2})\)~\cite{qian2025fundamental}. Second, all ensemble particles remain confined to the affine span of the initial ensemble~\cite{iglesias2013ensemble}. This subspace property can be beneficial as an implicit regularizer, but it also prevents recovery of directions not represented by the initial ensemble, which is restrictive when the ensemble size \(J\ll N_\theta\).

\subsection{Contributions}

Our main contributions are as follows:

\begin{enumerate}
\item 
We propose an inflated projected-dropout EKI method based on two modifications of standard EKI that preserve its parallel structure. Covariance inflation sustains ensemble spread, whereas projected dropout enables the ensemble mean to move beyond the initial affine span.

\item 
We provide a discrete-time, finite-ensemble analysis for linear inverse problems. The error in initial ensemble subspace decays exponentially at a rate independent of the condition number of the restricted Hessian, while the expected error in the complementary subspace also contracts exponentially. Meanwhile, the limiting covariance also provides a low-rank representation of posterior uncertainty.

\item Experiments on ill-conditioned linear systems, nonlinear benchmarks with curved or multimodal objective landscapes, and a PDE-constrained inverse problem for barotropic flow demonstrate accelerated convergence, sustained ensemble spread, and improved out-of-subspace recovery.
\end{enumerate}

\subsection{Literature Review}

\subsubsection{Ensemble Kalman inversion}

From an optimization perspective, ensemble Kalman
inversion~\cite{chen2012ensemble,emerick2013investigation,iglesias2013ensemble}
seeks to minimize a least-squares objective through repeated
ensemble Kalman updates.
Its convergence is studied through the long-time
behavior, as \(t\to\infty\), of the associated continuous-time EKI dynamics~\cite{schillings2017analysis,bungert2023complete}.
These methods are particularly attractive for expensive black-box
forward models whose derivatives are difficult to obtain or
unreliable~\cite{lea2000sensitivity,wang2014least}.
Applications include geophysics and reservoir history
matching~\cite{gu2007iterative,hendricks2008real,aanonsen2009ensemble},
turbulent flows~\cite{gao2021bi,zhang2022ensemble,shin2024multi}, fluid-structure interaction~\cite{cao2022bayesian}, and climate models~\cite{dunbar2021calibration}. 
In fluid applications, ensemble averaging can also provide a smoothing effect that
mitigates sensitivity to small-scale
oscillations~\cite[Lemma 2]{UKI}\cite{dunbar2022ensemble}. 
Adaptive surrogate models have also been combined
with Kalman inversion to reduce the cost of repeated
forward evaluations while controlling surrogate
error~\cite{yan2019adaptive,gao2024adaptive}.
The present work focuses on the deterministic square-root implementations~\cite{tippett2003ensemble}, 
including ensemble adjustment Kalman inversion
(EAKI) and ensemble transform Kalman inversion
(ETKI)~\cite{huang2022efficient,chada2022convergence}, based on the ensemble adjustment and ensemble transform Kalman filters~\cite{anderson2001ensemble,bishop2001adaptive},
respectively. 
Both reproduce the Kalman analysis covariance exactly, which facilitates the discrete-time, finite-ensemble analysis developed here.

From a sampling perspective, Kalman inversion methods have also been developed for Bayesian inverse problems. In that setting, the artificial time is often normalized to \([0,1]\), so that each update incorporates only a fraction of the likelihood. This gives rise to a homotopy or transport interpretation, in which the prior is deformed into the posterior over one unit of artificial time~\cite{emerick2013ensemble,ding2021ensemble,iglesias2021adaptive}. More generally, these methods are related to measure transport~\cite{reich2013nonparametric,law2015data,marzouk2017sampling,ruchi2019transform,pidstrigach2023affine}, sequential Monte Carlo methods~\cite{del2006sequential,doucet2001sequential}, and non-Gaussian extensions~\cite{baptista2024representation}. The present work instead considers a covariance-inflated artificial dynamical system whose long-time stationary behavior yields a Gaussian approximation of the posterior. This viewpoint is closely related to Fisher--Rao gradient flows under Gaussian~\cite{chopin2023connection,chen2023sampling}\cite[Remark 1]{huang2022efficient} or Gaussian mixture~\cite{chen2024efficient,che2025stable,che2026adaptive,yang2026polarized} approximations.

\subsubsection{Ensemble collapse and subspace restriction}
Standard EKI suffers from ensemble collapse and algebraically
slow convergence~\cite{schillings2017analysis,bungert2023complete,liu2025dropout}.
Existing approaches to acceleration include optimization-inspired
updates~\cite{kovachki2019ensemble,vernon2025nesterov},
covariance inflation, adaptive step sizes, and stabilized
continuous-time dynamics~\cite{chada2022convergence,weissmann2022gradient,armbruster2022stabilization,huang2022efficient}.
The covariance-inflated dynamics studied here follow this direction.

A separate limitation is that standard EKI remains confined to the
affine span of the initial ensemble~\cite{iglesias2013ensemble}.
Although this restriction provides implicit
regularization~\cite{parzer2022convergence}, it prevents recovery
of components outside that span.
Strategies for overcoming it include
localization~\cite{tong2023localized,al2024non}
and dropout perturbations~\cite{liu2025dropout}, the latter admitting
exponential convergence guarantees in expectation under suitable
assumptions.
Our projected dropout correction combines exploration beyond the
initial affine span with the covariance update of inflated EAKI.
It modifies only the ensemble mean, using observation responses
orthogonal to those represented by the original ensemble,
while preserving the analysis covariance for the current
predicted ensemble.

\subsection{Organization}

The remainder of the paper is organized as follows. In \cref{sec:inverse_problem}, we introduce the inverse problem, its Bayesian formulation, and the weighted nonlinear least-squares objective. In \cref{sec:inflated_EKI}, we introduce the methodology. In \cref{sec:theory}, we present a discrete-time, finite-ensemble convergence analysis in the linear setting. Numerical experiments are presented in \cref{sec:numerics}, followed by concluding remarks in \cref{sec:conclusion}.

\section{Inverse Problem and Weighted Nonlinear Least Squares}
\label{sec:inverse_problem}

We consider the inverse problem of recovering an unknown parameter
\(\theta \in \mathbb{R}^{N_\theta}\) from noisy observations
\(y \in \mathbb{R}^{N_y}\) related by~\cite{stuart2010inverse}
\begin{equation*}
\label{eq:KI}
    y = \mathcal{G}(\theta) + \eta.
\end{equation*}
Here \(\mathcal{G}\) denotes the forward model mapping parameters to observables, and \(\eta\) denotes observational noise. We assume that \(\eta\) is a realization of a mean-zero Gaussian random variable with covariance \(\Sigma_\eta\), namely,
$\eta \sim \mathcal{N}(0,\Sigma_\eta)$.
The corresponding data-misfit function is
\begin{equation}
\label{eq:data-misfit}
\Phi(\theta,y)
=
\frac12 (y-\mathcal{G}(\theta))^T \Sigma_\eta^{-1}(y-\mathcal{G}(\theta))
=
\frac12 \left\|\Sigma_\eta^{-1/2}(y-\mathcal{G}(\theta))\right\|^2.
\end{equation}
This is a weighted nonlinear least-squares objective, and minimizing \(\Phi\) yields an estimate of the unknown parameter \(\theta\).

Inverse problems are often ill posed and therefore require regularization. In the Bayesian framework, regularization is introduced by incorporating prior information on the unknown parameter. Specifically, \(\theta\) and \(y\) are treated as random variables. Given a prior density \(\rho_{\rm prior}(\theta)\), the inverse problem is formulated in terms of the posterior density \(\rho_{\rm post}(\theta)\) conditioned on the observation \(y\)~\cite{kaipio2006statistical,dashti2013bayesian,dashti2013map}:
\begin{equation*}
\label{eq:KI3}
\rho_{\rm post}(\theta)
=
\frac{1}{Z(y)}
e^{-\Phi(\theta,y)}\rho_{\rm prior}(\theta),
\end{equation*}
where \(Z(y)=
    \int e^{-\Phi(\theta,y)}\rho_{\rm prior}(\theta)\dd\theta\) is the normalization constant.
We focus on the case in which the prior \(\rho_{\rm prior}\) is Gaussian, or is approximated by a Gaussian distribution with mean \(r_0\) and covariance \(\Sigma_0\). In this case, the posterior density can be written as
\begin{equation}
\label{eq:posterior}
    \rho_{\rm post}(\theta)
    \propto e^{-\Phi_R(\theta,y)},
    \qquad
    \Phi_R(\theta,y)
    =
    \Phi(\theta,y)
    +
    \frac12\left\|\Sigma_0^{-1/2}(\theta-r_0)\right\|^2.
\end{equation}
The function \(\Phi_R\) is the Tikhonov-regularized version of the data misfit \(\Phi\)~\cite{chada2020tikhonov}, and minimizing \(\Phi_R\) yields the maximum a posteriori estimate.
This observation shows that both the unregularized inverse problem \eqref{eq:data-misfit} and the regularized Bayesian inverse problem \eqref{eq:posterior} can be treated within a unified weighted nonlinear least-squares framework. Accordingly, in the remainder of the paper we work with the general formulation
\begin{equation}
\label{eq:nlsq}
    \min_{\theta \in \R^{N_\theta}}
    \frac12 (x-\mathcal{F}(\theta))^T\Sigma_\nu^{-1}(x-\mathcal{F}(\theta))
    =
    \min_{\theta\in \R^{N_\theta}}\;\frac12 \left\|\Sigma_\nu^{-1/2}(x-\mathcal{F}(\theta))\right\|^2,
\end{equation}
where \((x,\mathcal{F},\Sigma_\nu)\) denotes either the original inverse problem or its augmented regularized formulation.

\section{Inflated Ensemble Kalman Inversion with Projected Dropout}
\label{sec:inflated_EKI}

In this section, we present our proposed ensemble Kalman inversion method, which combines covariance inflation with a projected dropout mean correction. We first introduce covariance-inflated mean-field dynamics and derive the corresponding ensemble adjustment Kalman inversion (EAKI) scheme. We then develop the projected dropout correction, which enables the ensemble mean to explore beyond the affine span of the initial ensemble without changing the covariance produced by the EAKI analysis step.

\subsection{Inflated Mean-Field Dynamics and EAKI}
\label{sec:EKI}
Ensemble Kalman inversion addresses the inverse problem~\eqref{eq:nlsq}
by introducing an artificial dynamical system and repeatedly assimilating
the same observed data using ensemble Kalman filtering.
To mitigate ensemble collapse and improve exploration within the ensemble
subspace, we replace the identity evolution in \eqref{eq:intro_dynamics}
with the following covariance-inflated mean-field dynamics:
\begin{subequations}
\label{eq:inflation_dynamics}
\begin{align}
&\textrm{evolution:}    &&\theta_{n+1} = m_{n} + \sqrt{\frac{1}{1 - \Dt}} (\theta_{n} - \mean_n), \label{eq:evolve} \\
  &\textrm{observation:}  &&\by_{n+1} = \bG(\theta_{n+1}) + \nu_{n+1}. \label{eq:obser}
\end{align}
\end{subequations}
Here \(\Delta\tau\in(0,1)\) denotes the artificial time step.
Let \(Y_n:=(\by_1^\dagger,\by_2^\dagger,\dots,\by_n^\dagger)\)
denote the observation history, with \(\by_k^\dagger=x\) for every \(k\).
The conditional law of \(\theta_n\) given \(Y_n\) is the filtering
distribution at iteration \(n\), and
\(m_n=\mathbb E[\theta_n\mid Y_n]\) denotes its mean.
Equation~\eqref{eq:evolve} is therefore of mean-field type, since it
depends on the filtering distribution through its mean.
The artificial observation noises are independent Gaussian random
variables with
\begin{equation}
\label{eq:hyperparameters}
\nu_{n+1}\sim\mathcal N(0,\Sigma_{\nu,n+1}), \qquad
\Sigma_{\nu,n+1}=\frac{1}{\Dt}\Sigma_\nu.
\end{equation}
The scaling with respect to \(\Delta\tau\) is crucial for ensuring that the resulting filtering distribution provides a consistent Gaussian approximation of the posterior uncertainty; see \cref{thm:EAKI_ETKI_convergence}.

We approximate the filtering distribution of \(\theta_n\) given
\(Y_n\) by an ensemble of \(J\) particles \(\{\theta_n^j\}_{j=1}^J\). At the ensemble level, we use \(m_n\) to denote the empirical
mean \(J^{-1}\sum_{j=1}^J\theta_n^j\), which approximates the
conditional mean \(\mathbb E[\theta_n\mid Y_n]\).
Applying the inflated evolution~\eqref{eq:evolve} to each particle
yields a predicted ensemble that approximates the conditional
distribution of \(\theta_{n+1}\) given \(Y_n\):
\begin{equation}
\label{eq:inflated_prediction}
    \widehat\theta_{n+1}^{j}
    =
    m_n+
    \sqrt{\frac{1}{1-\Dt}}
    (\theta_n^j-m_n),
    \qquad j=1,\dots,J ,
\end{equation}
with empirical mean $\widehat{m}_{n+1}$ and covariance $\widehat{C}_{n+1}$.

In the analysis step, we consider the predicted observation
\(x_{n+1}=\F(\theta_{n+1})+\nu_{n+1}\) defined
by~\eqref{eq:obser}. We approximate the conditional joint distribution
of \((\theta_{n+1},x_{n+1})\) given \(Y_n\) by the Gaussian distribution
\begin{equation}
\label{eq:KF_joint}
\mathcal N\!\left(
\begin{bmatrix}
\widehat m_{n+1}\\[2mm]
\widehat x_{n+1}
\end{bmatrix},
\begin{bmatrix}
\widehat C_{n+1} & \widehat C_{n+1}^{\theta x}\\
(\widehat C_{n+1}^{\theta x})^T & \widehat C_{n+1}^{xx}
\end{bmatrix}
\right),
\end{equation}
where the predictive moments satisfy
\begin{equation}
\label{eq:KF_joint2}
\begin{split}
    \widehat x_{n+1}
    &= \mathbb E[x_{n+1}\mid Y_n]
     = \mathbb E[\F(\theta_{n+1})\mid Y_n], \\
    \widehat C_{n+1}^{\theta x}
    &= \operatorname{Cov}[\theta_{n+1},x_{n+1}\mid Y_n]
     = \operatorname{Cov}[\theta_{n+1},\F(\theta_{n+1})\mid Y_n], \\
    \widehat C_{n+1}^{xx}
    &= \operatorname{Cov}[x_{n+1}\mid Y_n]
     = \operatorname{Cov}[\F(\theta_{n+1})\mid Y_n] + \Sigma_{\nu,n+1}.
\end{split}
\end{equation}
Conditioning the Gaussian approximation~\eqref{eq:KF_joint} on
\(x_{n+1}=x\) yields a Gaussian approximation to the filtering
distribution of \(\theta_{n+1}\) given \(Y_{n+1}\).
The mean and covariance of this approximation are given by the
Kalman update:
\begin{subequations}
\label{eq:sds}
\begin{align}
    m_{n+1}
    &= \widehat m_{n+1}
    + \widehat C_{n+1}^{\theta x}
      (\widehat C_{n+1}^{xx})^{-1}
      (x-\widehat x_{n+1}),
    \label{eq:KF_analysis_mean}
    \\
    C_{n+1}
    &= \widehat C_{n+1}
    - \widehat C_{n+1}^{\theta x}
      (\widehat C_{n+1}^{xx})^{-1}
      (\widehat C_{n+1}^{\theta x})^T.
    \label{eq:KF_analysis_cov}
\end{align}
\end{subequations}

In ensemble Kalman methods, the predictive moments
in \eqref{eq:KF_joint2} are estimated empirically, and the filtering
distribution of \(\theta_{n+1}\) given \(Y_{n+1}\) is approximated
by an updated ensemble.
For the predicted ensemble~\eqref{eq:inflated_prediction}, we evaluate
the forward map and compute the empirical means:
\begin{equation}
\label{eq:EKI-Fx-analysis}
\widehat x_{n+1}^j=
\F(\widehat\theta_{n+1}^j),
\quad
\widehat x_{n+1}
=
\frac1J\sum_{j=1}^J\widehat x_{n+1}^j .
\end{equation}
We then form the normalized ensemble anomaly matrices
\begin{subequations}
\label{eq:EKF_square_root}
\begin{align}
    \widehat Z_{n+1}
    &=
    \frac{1}{\sqrt{J-1}}
    \begin{bmatrix}
        \widehat\theta_{n+1}^1-m_{n}
        &
        \cdots
        &
        \widehat\theta_{n+1}^J-m_{n}
    \end{bmatrix} \in\R^{N_\theta \times J},
    \label{eq:hat_Z}\\
    \widehat Y_{n+1}
    &=
    \frac{1}{\sqrt{J-1}}
    \begin{bmatrix}
        \widehat x_{n+1}^1-\widehat x_{n+1}
        &
        \cdots
        &
        \widehat x_{n+1}^J-\widehat x_{n+1}
    \end{bmatrix} \in\R^{N_x \times J}.\label{eq:hat_Y}
\end{align}
\end{subequations}
Here we used $\widehat{m}_{n+1} = m_n$. These matrices yield the empirical covariance estimates in~\eqref{eq:KF_joint2}
\begin{equation}
\label{eq:EKF_covariances}
    \widehat C_{n+1}
    =
    \widehat Z_{n+1}\widehat Z_{n+1}^T,
    \qquad
    \widehat C_{n+1}^{\theta x}
    =
    \widehat Z_{n+1}\widehat Y_{n+1}^T,
    \qquad
    \widehat C_{n+1}^{xx}
    =
    \widehat Y_{n+1}\widehat Y_{n+1}^T+\Sigma_{\nu,n+1},
\end{equation}
where the observation noise covariance is included analytically.

The analysis step constructs an updated ensemble
\(\{\theta_{n+1}^j\}_{j=1}^J\) to approximate the filtering distribution.
A key computational advantage of ensemble Kalman methods is that
the particle updates can be implemented using ensemble anomaly matrices,
without explicitly forming or storing the full
\(N_\theta\times N_\theta\) parameter covariance matrix.
This is particularly useful in high-dimensional problems
with \(J\ll N_\theta\).
We focus on ensemble adjustment Kalman inversion (EAKI), which is based
on the ensemble adjustment Kalman filter~\cite{anderson2001ensemble}.
EAKI updates the mean according to \eqref{eq:KF_analysis_mean} and
transforms the ensemble anomalies to reproduce the Kalman analysis
covariance~\eqref{eq:KF_analysis_cov}.
Using the Woodbury identity, this covariance can be written as
\begin{equation}
\label{eq:analysis_cov_square_root_form}
    C_{n+1}
    =
    \widehat Z_{n+1}
    \Bigl(
    I+\widehat Y_{n+1}^T\Sigma_{\nu,n+1}^{-1}\widehat Y_{n+1}
    \Bigr)^{-1}
    \widehat Z_{n+1}^T .
\end{equation}
EAKI realizes this covariance through a linear transformation acting
on the ensemble anomalies from the left. 
To construct this transformation, consider the compact singular value
decomposition of \(\widehat Z_{n+1}\) and the associated spectral
decomposition:
\begin{equation}
\label{eq:eaki-svd}
\begin{split}
   \text{Compact SVD:} \qquad
   &\widehat Z_{n+1}
   =
   U_Z\widehat S_Z^{1/2}V_Z^T,
   \\
   \text{Spectral decomposition:} \qquad
   &V_Z^T
   \Bigl(
    I+\widehat Y_{n+1}^T\Sigma_{\nu,n+1}^{-1}\widehat Y_{n+1}
   \Bigr)^{-1}
   V_Z
   =
   W_Z S_A W_Z^T.
\end{split}
\end{equation}
Here \(\widehat S_Z\) and \(S_A\) are nonsingular diagonal matrices of size
\(\operatorname{rank}(\widehat Z_{n+1})\). Define the transformation matrix 
\begin{equation}
\label{eq:EAKI_transform}
    A_{n+1}
    =
    U_Z\widehat S_Z^{1/2}
    W_Z S_A^{1/2}
    \widehat S_Z^{-1/2}U_Z^T.
\end{equation}
The EAKI analysis step is then
\begin{subequations}
\label{eq:EAKI-analysis}
\begin{align}
m_{n+1}
&=
m_{n}
+
\widehat C_{n+1}^{\theta x}
\left(\widehat C_{n+1}^{xx}\right)^{-1}
(x-\widehat x_{n+1}),
\label{eq:EAKI-mean-analysis}
\\
\theta_{n+1}^{j}
&=
m_{n+1}
+
A_{n+1}
(\widehat\theta_{n+1}^{j}-\widehat m_{n+1}),
\qquad
j=1,\dots,J .
\label{eq:EAKI-particle-analysis}
\end{align}
\end{subequations}
By construction, the normalized anomaly matrix \(Z_{n+1}\) of the
updated ensemble satisfies
\begin{equation}
\label{eq:cov-eaki}
    Z_{n+1}Z_{n+1}^T
    =
    \widehat Z_{n+1}
    \Bigl(
    I+\widehat Y_{n+1}^T\Sigma_{\nu,n+1}^{-1}\widehat Y_{n+1}
    \Bigr)^{-1}
    \widehat Z_{n+1}^T
    =
    \widehat C_{n+1}
    -
    \widehat C_{n+1}^{\theta x}
    \bigl(\widehat C_{n+1}^{xx}\bigr)^{-1}
    \bigl(\widehat C_{n+1}^{\theta x}\bigr)^T .
\end{equation}
Thus EAKI is a deterministic ensemble implementation of the Gaussian conditioning update \eqref{eq:sds}.

\paragraph{Subspace property.}
Both standard EAKI and its covariance-inflated variant remain confined
to the affine span of the initial ensemble.
Let \(\mathcal Q:=\Range(Z_0)\).
For both methods, the prediction step~\eqref{eq:inflated_prediction} preserves the anomaly subspace: $\Range(\widehat Z_{n+1})=\Range(Z_n)$,
since it either leaves the anomalies unchanged or rescales them by
a nonzero factor.
Moreover, the EAKI mean increment~\eqref{eq:EAKI-mean-analysis}
lies in \(\Range(\widehat Z_{n+1})\).
By \eqref{eq:EAKI_transform}, \(A_{n+1}\) has the same column space
as \(\widehat Z_{n+1}\), and its restriction to this subspace is
invertible. Consequently,
$\Range(Z_{n+1})
=
\Range(A_{n+1}\widehat Z_{n+1})
=
\Range(\widehat Z_{n+1})$.
It follows by induction that
\(\theta_n^j\in m_0+\mathcal Q\) and
\(\Range(Z_n)=\mathcal Q\) for all \(n\geq0\).
Inflation therefore increases ensemble spread without enlarging
the affine span of the initial ensemble.

\subsection{Projected Dropout Mean Correction}
\label{ssec:dropout_method}

Covariance inflation expands the predicted anomalies within \(\mathcal Q\), but the ensemble mean remains in \(m_0+\mathcal Q\). 
To allow the mean to  explore beyond this  affine space, we augment its update with directions generated by coordinate-wise dropout~\cite{liu2025dropout}.
Setting selected coordinates of an anomaly~\eqref{eq:hat_Z} to zero can change its direction and move it outside \(\mathcal Q\). We evaluate these
candidate directions through the forward map and retain the parts of their observation responses that are not represented by the original ensemble.

Specifically, at iteration \(n\), the diagonal mask \(R_n\in\R^{N_\theta\times N_\theta}\) is drawn independently of earlier masks from independent \(\operatorname{Bernoulli}(\lambda)\) entries, where \(0<\lambda<1\), with the entire mask resampled if it
is zero. Applying the same mask to all predicted anomalies~\eqref{eq:hat_Z} preserves their zero sum and gives the centered perturbations
\(\widetilde\tau_{n+1}^j=R_n(\widehat\theta_{n+1}^j-m_n)\).
Their forward responses
\(\widetilde x_{n+1}^j=\F(m_n+\widetilde\tau_{n+1}^j)\) 
have mean
\(\widetilde x_{n+1}=J^{-1}\sum_{j=1}^J\widetilde x_{n+1}^j\).
The normalized parameter and observation anomaly matrices are then
\begin{subequations}
\label{eq:dropout-masked-anomalies}
\begin{align}
\widetilde Z_{n+1}
&=R_n\widehat Z_{n+1}
=\frac{1}{\sqrt{J-1}}
\begin{bmatrix}
\widetilde\tau_{n+1}^1 , \cdots , \widetilde\tau_{n+1}^J
\end{bmatrix},
\\
\widetilde Y_{n+1}
&=\frac{1}{\sqrt{J-1}}
\begin{bmatrix}
\widetilde x_{n+1}^1-\widetilde x_{n+1} , \cdots ,
\widetilde x_{n+1}^J-\widetilde x_{n+1}
\end{bmatrix}.
\end{align}
\end{subequations}

The columns of \(\widehat Y_{n+1}\) in~\eqref{eq:hat_Y} and $\widetilde Y_{n+1}$ describe the observation changes sampled by the original and perturbed ensembles, respectively. To isolate complementary observation responses, we project \(\widetilde Y_{n+1}\) onto the orthogonal complement of \(\Range(\widehat Y_{n+1})\) in the observation-noise metric. This metric agrees with that of the least-squares objective, since 
\(\Sigma_{\nu,n+1}^{-1}=\Dt\Sigma_\nu^{-1}\). The projected observation anomalies and the associated parameter directions are
\begin{subequations}
\label{eq:dropout-projected-residuals}
\begin{align}
L_{n+1}
&=\left(
\widehat Y_{n+1}^T\Sigma_{\nu,n+1}^{-1}\widehat Y_{n+1}
\right)^\dagger
\widehat Y_{n+1}^T\Sigma_{\nu,n+1}^{-1}\widetilde Y_{n+1},
\\
\widetilde Y_{n+1}^{\perp}
&=\widetilde Y_{n+1}-\widehat Y_{n+1}L_{n+1},
\\
\widetilde Z_{n+1}^{\perp}
&=\widetilde Z_{n+1}-\widehat Z_{n+1}L_{n+1}.
\end{align}
\end{subequations}
Here \((\cdot)^\dagger\) denotes the Moore--Penrose pseudoinverse.
The common coefficients \(L_{n+1}\) keep the parameter directions $\widetilde Z_{n+1}^{\perp}$
paired with the observation anomalies $\widetilde Y_{n+1}^{\perp}$. Moreover, after whitening by \(\Sigma_{\nu,n+1}^{-1/2}\), the original and projected observation blocks are orthogonal: 
\begin{equation}
\label{eq:dropout-observation-orthogonality}
\widehat Y_{n+1}^T\Sigma_{\nu,n+1}^{-1}
\widetilde Y_{n+1}^{\perp}=0.
\end{equation}

We construct the augmented anomaly matrices by combining the original and projected blocks:
\begin{equation}
\label{eq:joint-augmented-anomalies}
\begin{aligned}
Z_{n+1}^{\rm aug}
=
\left[
\widehat Z_{n+1},
\widetilde Z_{n+1}^{\perp}
\right], \quad
Y_{n+1}^{\rm aug}
=
\left[
\widehat Y_{n+1},
\widetilde Y_{n+1}^{\perp}
\right].
\end{aligned}
\end{equation}
In whitened observation coordinates, orthogonality~\eqref{eq:dropout-observation-orthogonality} separates the two contributions to the joint Kalman gain
\begin{align}
\label{eq:augmented-Kalman-gain}
K_{n+1}^{\rm aug}
:=&Z_{n+1}^{\rm aug}
(Y_{n+1}^{\rm aug})^T
\left(
Y_{n+1}^{\rm aug}(Y_{n+1}^{\rm aug})^T
+\Sigma_{\nu,n+1}
\right)^{-1} \\
=&
\underbrace{
\widehat Z_{n+1}\widehat Y_{n+1}^T
\left(
\widehat Y_{n+1}\widehat Y_{n+1}^T
+\Sigma_{\nu,n+1}
\right)^{-1}
}_{ \rm original} 
+
\underbrace{
\widetilde Z_{n+1}^{\perp}
(\widetilde Y_{n+1}^{\perp})^T
\left(
\widetilde Y_{n+1}^{\perp}
(\widetilde Y_{n+1}^{\perp})^T
+\Sigma_{\nu,n+1}
\right)^{-1}
}_{\rm projected}. \notag
\end{align}
The original block
responds to residual components in its observation span,
while the projected block responds to components in the
complementary span sampled by the projected anomalies.
Applying this gain to the residual at the current mean $m_n$ gives
\begin{equation}
\label{eq:dropout-mean-correction}
m_{n+1}
=
m_n
+
Z_{n+1}^{\rm aug}
(Y_{n+1}^{\rm aug})^T
\left(
Y_{n+1}^{\rm aug}(Y_{n+1}^{\rm aug})^T
+\Sigma_{\nu,n+1}
\right)^{-1}
(x-\F(m_n)).
\end{equation}
For a linear forward map, \(\F(m_n)=\widehat x_{n+1}\), so this residual coincides with that used in the ordinary EAKI mean update~\eqref{eq:EAKI-mean-analysis}.

Finally, for the ensemble particle update, we use the augmented
blocks only to update the mean, allowing it to move beyond
the initial affine span.
The anomaly transform \(A_{n+1}\) is computed from the original
anomalies via~\eqref{eq:EAKI_transform}, preserving the EAKI
analysis covariance for the current predicted ensemble.
The resulting particle update is
\begin{equation}
\label{eq:dropout-particle-translation}
\theta_{n+1}^j
=m_{n+1}
+A_{n+1}\left(
\widehat\theta_{n+1}^j-m_n
\right),
\quad j=1,\ldots,J.
\end{equation}
The linear theory in \cref{ssec:dropout_theory} shows that, under the coordinate-coverage
condition~\eqref{eq:dropout-coordinate-coverage}, this mean correction
preserves EAKI's exponential convergence within the initial
ensemble subspace $\mathcal Q$ and yields exponential mean-square decay
of the complementary observable error.

\begin{algorithm}[H]
\caption{Inflated projected-dropout EAKI}
\label{alg:projected_dropout_mean_correction}
\begin{algorithmic}[1]
\State \textbf{Input:} initial ensemble \(\{\theta_0^j\}_{j=1}^J\), data \(x\), forward map \(\F\), covariance \(\Sigma_\nu\), step size \(\Dt\), dropout rate \(1-\lambda\), and iteration count \(N\).
\For{\(n=0,1,\ldots,N-1\)}
    \State Apply the inflated prediction \eqref{eq:inflated_prediction} and form \(\widehat Z_{n+1}\) and \(\widehat Y_{n+1}\).
    \State Compute \(A_{n+1}\) from the original predicted anomalies using \eqref{eq:eaki-svd} and  \eqref{eq:EAKI_transform}.
    \State Generate \(R_n\), resampling if the resulting mask is all zeros, and compute the masked anomalies via \eqref{eq:dropout-masked-anomalies}.
    \State Compute the projected and augmented anomalies by \eqref{eq:dropout-projected-residuals} and \eqref{eq:joint-augmented-anomalies}.
    \State Update the mean by \eqref{eq:dropout-mean-correction} and the particles by \eqref{eq:dropout-particle-translation}.
\EndFor
\State \textbf{Output:} corrected ensemble \(\{\theta_N^j\}_{j=1}^J\) and mean \(m_N\).
\end{algorithmic}
\end{algorithm}

\begin{remark}[Geometry of the projection]
The observation projection~\eqref{eq:dropout-projected-residuals} measures parameter directions through their
noise-weighted effects on the forward response. If
\(\F(\theta)=F\theta\), then
\(\widehat Y_{n+1}=F\widehat Z_{n+1}\) and
\(\widetilde Y_{n+1}^{\perp}=F\widetilde Z_{n+1}^{\perp}\).
Since \(\Sigma_{\nu,n+1}^{-1}=\Dt\Sigma_\nu^{-1}\), the weighted
orthogonality above becomes
\[
\widehat Z_{n+1}^{T}F^T\Sigma_\nu^{-1}F
\widetilde Z_{n+1}^{\perp}=0.
\]
Here \(F^T\Sigma_\nu^{-1}F\) is the Hessian of the weighted
least-squares objective. Thus the projected parameter anomalies are
orthogonal to the original ones in the Hessian form, which accounts
for forward sensitivity and observation noise.

\end{remark}

The proposed method requires \(2J+1\) evaluations of the
forward map per iteration.
\Cref{alg:projected_dropout_mean_correction} summarizes
the complete algorithm.

\section{Theoretical Analysis}
\label{sec:theory}

In this section, we provide a discrete-time, finite-ensemble analysis
of the proposed EAKI method for linear inverse problems to clarify
the roles of covariance inflation and projected dropout.
We assume that the forward map takes the form
\begin{equation}
    \mathcal{F}(\theta)=F\theta,
\end{equation}
and the associated weighted least-squares objective is
$\Phi_R(m) = \frac12\left\|\Sigma_\nu^{-1/2}(x-Fm)\right\|^2$.
For an initial ensemble of size \(J\geq2\), let
\(m_0\in\mathbb{R}^{N_\theta}\),
\(C_0\in\mathbb{R}^{N_\theta\times N_\theta}\), and
\(Z_0\in\mathbb{R}^{N_\theta\times J}\)
denote its empirical mean, empirical covariance, and normalized
anomaly matrix, respectively, so that \(C_0=Z_0Z_0^T\).
Let \(r=\operatorname{rank}(C_0)\).
Since the ensemble anomalies are centered,
\(r\leq\min\{N_\theta,J-1\}\).
The covariance \(C_0\) admits the compact spectral decomposition
\begin{equation}
\label{eq:D0}
    C_0=QD_0Q^T,
\end{equation}
where \(Q\in\mathbb{R}^{N_\theta\times r}\) satisfies
\(Q^TQ=I_r\), and \(D_0\in\mathbb{R}^{r\times r}\) is diagonal
with strictly positive entries.
The columns of \(Q\) form an orthonormal basis of
the initial anomaly subspace.

In \cref{ssec:linear_theory}, we show that covariance inflation
yields exponential convergence of the ensemble mean within the
initial affine span, with a convergence factor independent of
conditioning.
In \cref{ssec:dropout_theory}, we analyze the full method and
show that projected dropout enables exploration beyond the initial
affine span and yields exponential convergence of the expected
least-squares objective value to its minimum, at a rate that may
depend on conditioning.

\subsection{Covariance inflation}
\label{ssec:linear_theory}

Without the dropout correction, the inflated EAKI updates
in~\cref{eq:inflated_prediction,eq:EAKI-analysis} remain confined
to the initial affine span \(m_0+\Range(Q)\).
To account for the initial covariance, we introduce the normalized
coordinates $\theta=m_0+QD_0^{1/2}z$.
In the \(z\)-coordinates, the initial ensemble has zero empirical
mean and identity empirical covariance.

The least-squares problem restricted to this affine space need not
have a unique minimizer, since \(\Range(Q)\) may contain nonzero
directions in \(\ker(F)\).
To separate observable and unobservable directions, we consider
the compact spectral decomposition of the least-squares Hessian
in the \(z\)-coordinates:
\begin{equation}
\label{eq:P_1_Gamma_P1}
D_0^{1/2}Q^TF^T\Sigma_\nu^{-1}FQD_0^{1/2}
=
P_1\Gamma P_1^T,
\end{equation}
where \(r'=\operatorname{rank}(FQ)\) and  $P_1\in\mathbb R^{r\times r'}$,    $P_1^TP_1=I_{r'}$,  $\Gamma=\operatorname{diag}(\gamma_1,\dots,\gamma_{r'})$,
$\gamma_i>0$.
Let \(P_2\in\mathbb R^{r\times(r-r')}\) have orthonormal columns
spanning the null space of this Hessian, and define the orthogonal
matrix $P=[P_1\;P_2]\in\mathbb R^{r\times r}$.
This decomposition separates the normalized coordinates into
observable and unobservable components, in the sense that
\begin{equation}
\label{eq:linear-null-directions}
    FQD_0^{1/2}P_2=0.
\end{equation}
Thus, writing $m=m_0+QD_0^{1/2}(P_1b+P_2c)$,  with  $b\in\mathbb R^{r'}$, $c\in\mathbb R^{r-r'}$,
the data misfit depends only on \(b\).
The full and reduced problems therefore have the same minimum value:
 \begin{equation*}
 \min_{z\in\mathbb R^{r}}
    \frac12
    \left\|
    \Sigma_\nu^{-1/2}
    \bigl(x-Fm_0-FQD_0^{1/2}z\bigr)
    \right\|^2 = 
    \min_{b\in\mathbb R^{r'}}
    \frac12
    \left\|
    \Sigma_\nu^{-1/2}
    \bigl(x-Fm_0-FQD_0^{1/2}P_1b\bigr)
    \right\|^2.
    \end{equation*}
The reduced problem in \(b\) has Hessian \(\Gamma\succ0\)
and hence the unique minimizer
\begin{equation}
\label{eq:b_star}
    b_\star
    =
    \Gamma^{-1}P_1^TD_0^{1/2}Q^TF^T
    \Sigma_\nu^{-1}(x-Fm_0).
\end{equation}
We select the minimum-norm representative \(z_\star=P_1b_\star\), and define the distinguished minimizer
over the initial affine span by $m_\star=m_0+QD_0^{1/2}P_1b_\star$.


\begin{theorem}[Inflated EAKI in the linear case]
\label{thm:EAKI_ETKI_convergence}
Assume that \(0<\Dt<1\), \(\Sigma_\nu\succ0\), and \(\F(\theta)=F\theta\). Then, for the inflated EAKI, the following statements hold:
\begin{enumerate}
    \item \textbf{Mean and covariance recursions.} The empirical mean and covariance satisfy
    \begin{subequations}
    \label{eq:linear-moment-recursions}
    \begin{align}
        m_{n+1}
        &=m_n+C_nF^T
          \left(
              FC_nF^T+\frac{1-\Dt}{\Dt}\Sigma_\nu
          \right)^{-1}
          (x-Fm_n),
        \label{eq:linear-mean-recursion}
        \\
        C_{n+1}
        &=\frac1{1-\Dt}
          \left[
              C_n-C_nF^T
              \left(
                  FC_nF^T+\frac{1-\Dt}{\Dt}\Sigma_\nu
              \right)^{-1}
              FC_n
          \right].
        \label{eq:linear-cov-recursion}
    \end{align}
    \end{subequations}
    
    \item \textbf{Exponential covariance convergence in observable directions.} 
    For every \(n\ge0\), $C_n=QD_nQ^T$  with \(D_n\succ0\) satisfying
    \begin{equation}
    \label{eq:linear-precision-closed}
        D_n^{-1}=Q^T F^T\Sigma_{\nu}^{-1} FQ +(1-\Dt)^n\left(D_0^{-1}-Q^T F^T\Sigma_{\nu}^{-1} FQ\right).
    \end{equation}
Hence the pseudoinverse \(C_n^\dagger=QD_n^{-1}Q^T\)
converges exponentially to the projected Hessian
\(QQ^TF^T\Sigma_\nu^{-1}FQQ^T\).
Moreover,
    \begin{equation}
    \label{eq:linear-covariance-closed}
        D_n
        =D_0^{1/2}
          \left(
              P_1\Lambda_nP_1^T
              +(1-\Dt)^{-n}P_2P_2^T
          \right)
          D_0^{1/2},
    \end{equation}
    where $\Lambda_n :=\left[(1-\Dt)^nI_{r'}+(1-(1-\Dt)^n)\Gamma\right]^{-1}.$

    \item \textbf{Exponential mean convergence in the initial ensemble subspace.}
    \begin{equation}
    \label{eq:linear-mean-representation}
        m_n=m_0+ Q D_0^{1/2} P_1 b_n,
        \quad
        b_n-b_\star=-(1-\Dt)^n\Lambda_nb_\star,
    \end{equation}
    and
    \begin{equation}
    \label{eq:linear-mean-rate}
        \|b_n-b_\star\|_2
        \le
        \frac{\|b_\star\|_2}
             {\min\{1,\lambda_{\min}(\Gamma)\}}
        (1-\Dt)^n.
    \end{equation}
Thus the geometric convergence factor is
\(1-\Dt\), independently of the condition number of \(\Gamma\),
although the multiplicative constants may depend on \(\Gamma\).
\end{enumerate}
\end{theorem}

\begin{proof}
For $(i)$, the prediction step~\eqref{eq:inflated_prediction} gives, in the linear case,
\[
    \widehat m_{n+1}=m_n,
    \quad
    \widehat x_{n+1}=Fm_n,
    \quad
    \widehat Z_{n+1}=(1-\Dt)^{-1/2}Z_n,
    \quad
    \widehat Y_{n+1}=F\widehat Z_{n+1}.
\]
Substituting these identities and the noise scaling
\eqref{eq:hyperparameters} into the EAKI mean and
covariance updates \eqref{eq:EAKI-mean-analysis} and \eqref{eq:cov-eaki} gives
\eqref{eq:linear-moment-recursions}. Since \(\operatorname{rank}(C_n)\leq J-1\), the covariance
is singular whenever \(J\leq N_\theta\).

For $(ii)$, the base case follows from \eqref{eq:D0}.
Suppose \(C_n=QD_nQ^T\) with \(D_n\succ0\). Applying the Woodbury identity to
the covariance update \eqref{eq:linear-cov-recursion} gives
\begin{align*}
    C_{n+1}
      &=Q\left((1-\Dt) D_n^{-1}+\Dt Q^T F^T\Sigma_{\nu}^{-1} FQ\right)^{-1}Q^T.
\end{align*}
Thus $C_{n+1}=QD_{n+1}Q^T$ with $D_{n+1}^{-1}=(1-\Dt) D_n^{-1}+\Dt Q^T F^T\Sigma_{\nu}^{-1} FQ \succ 0$. Solving the recursion gives
\eqref{eq:linear-precision-closed}. Multiplying that identity on the left and right by \(D_0^{1/2}\)
and using \eqref{eq:P_1_Gamma_P1}, we obtain
\begin{align*}
    D_0^{1/2}D_n^{-1}D_0^{1/2}
    &=(1-\Dt)^nI_r+(1-(1-\Dt)^n)P_1\Gamma P_1^T =P_1\Lambda_n^{-1}P_1^T+(1-\Dt)^nP_2P_2^T.
\end{align*}
Inverting this identity yields \eqref{eq:linear-covariance-closed}. Since \(C_n^\dagger=QD_n^{-1}Q^T\), its exponential convergence
follows directly from \eqref{eq:linear-precision-closed}.

For $(iii)$, we rewrite the increment in the mean update
\eqref{eq:linear-mean-recursion} using the following identity:
\begin{equation}
\label{eq:linear-gain-identity}
\begin{split}
C_nF^T
\left(FC_nF^T+\frac{1-\Dt}{\Dt}\Sigma_\nu\right)^{-1} 
      &=  \Dt C_{n+1}F^T\Sigma_\nu^{-1}
      = \Dt Q D_{n+1}Q^T F^T\Sigma_\nu^{-1} \\
      &=  \Dt Q D_0^{1/2} P_1\Lambda_{n+1} P_1^T D_0^{1/2} Q^TF^T \Sigma_\nu^{-1}.
\end{split}
\end{equation}
Here the first equality uses the covariance update~\eqref{eq:linear-cov-recursion}, the second uses $C_{n+1}=QD_{n+1}Q^T$, and the last follows from \eqref{eq:linear-covariance-closed} and the null-space identity~\eqref{eq:linear-null-directions} of $P_2$. 
Thus every mean increment lies in \(\Range(QD_0^{1/2}P_1)\),
so \(m_n=m_0+Q D_0^{1/2} P_1  b_n\) with \(b_0=0\).
Substituting this representation into~\eqref{eq:linear-mean-recursion} and using the definitions of
\(\Gamma\) and \(b_\star\) yields
\begin{equation}
    \label{eq:linear-b-recursion}
    \begin{split}
    b_{n+1}
    &=b_n+\Dt\Lambda_{n+1}
      P_1^T D_0^{1/2} Q^T F^T \Sigma_\nu^{-1} (x-Fm_0- FQ D_0^{1/2} P_1  b_n) \\
    &=b_n+\Dt\Lambda_{n+1}\Gamma(b_\star-b_n) \\
    &=b_{\star} + (1-\Dt)\Lambda_{n+1}\Lambda_n^{-1}(b_n-b_\star)
    \end{split}
\end{equation}
The last equality uses the recursion  
$\Lambda_{n+1}^{-1}
    =(1-\Dt)\Lambda_n^{-1}+\Dt\Gamma$,
 from part (ii), which implies 
$
    I-\Dt\Lambda_{n+1}\Gamma
    =(1-\Dt)\Lambda_{n+1}\Lambda_n^{-1}.
$ 
Since \(\Lambda_0=I_{r'}\) and \(b_0=0\), iteration gives
\[
    b_n-b_\star=-(1-\Dt)^n\Lambda_nb_\star,
\]
proving \eqref{eq:linear-mean-representation}.
Finally, every eigenvalue of
\(\Lambda_n^{-1}=(1-\Dt)^nI_{r'}+(1-(1-\Dt)^n)\Gamma\) is bounded below by
\(\min\{1,\lambda_{\min}(\Gamma)\}\).
Combining this bound with \eqref{eq:linear-mean-representation}
yields \eqref{eq:linear-mean-rate}.
\end{proof}

\begin{remark}[Unobservable covariance directions]
If the initial ensemble subspace contains nonzero directions
in the null space of \(F\), then \(r'<r\) in
\eqref{eq:P_1_Gamma_P1}.
The second term in the covariance formula
\eqref{eq:linear-covariance-closed} grows like \((1-\Dt)^{-n}\).
Consequently, the uncertainty represented by the ensemble in
these directions tends to infinity.
This behavior emphasizes the importance of incorporating the
prior distribution through the augmented formulation
\eqref{eq:posterior}, which regularizes these otherwise
unobservable directions.
\end{remark}

\subsection{Projected Dropout Correction}
\label{ssec:dropout_theory}
We next analyze the projected dropout mean
correction introduced in \cref{ssec:dropout_method}. 

\begin{theorem}[Inflated Projected Dropout EAKI in the linear case]
\label{thm:DropoutEAKI_ETKI}
Assume the linear setting and notation of
\cref{thm:EAKI_ETKI_convergence}. For each \(n\), independently draw a vector
with i.i.d. \(\operatorname{Bernoulli}(\lambda)\) entries, where
\(0<\lambda<1\), resampling the entire vector whenever it is zero, and let
\(R_n\) be the diagonal matrix formed from the accepted vector. Recall that the columns of \(Q\) form an orthonormal basis
of the initial anomaly subspace. Assume \(r'={\rm rank}(FQ)\geq1\), and define
\[
    \mathcal S:=\Range(\Sigma_\nu^{-1/2}FQ),
    \qquad
    \mathcal R:=\mathcal S^\perp\cap\Range(\Sigma_\nu^{-1/2}F),
\]
where \(\mathcal S^\perp\) is the orthogonal complement of $\mathcal S$ in $\R^{N_x}$. Let \(P_\perp\) denote the Euclidean projection onto \(\mathcal S^\perp\).
Suppose that the initial ensemble subspace satisfies the coordinate-coverage
condition
\begin{equation}
\label{eq:dropout-coordinate-coverage}
    \ell_Q
    :=
    \min_{1\leq i\leq N_\theta}e_i^TQQ^Te_i
    >0.
\end{equation}
Let \(m^*\) be a minimizer of
\[
    \min_{m\in\R^{N_\theta}}\Phi_R(m),
    \qquad
    \Phi_R(m)=\frac12\left\|\Sigma_\nu^{-1/2}(x-Fm)\right\|^2,
\]
and let \(b_\star\) be defined by~\eqref{eq:b_star}. Then, for every
\(n\geq0\), the following statements hold. Throughout, all expectations are taken only
over the dropout masks.

\begin{enumerate}
\item \textbf{Covariance preservation.}
The projected dropout mean correction leaves the covariance update unchanged.
Consequently, all identities and convergence estimates for the covariance in
\cref{thm:EAKI_ETKI_convergence} remain valid.

\item \textbf{In-subspace mean convergence.}
There exist \(b_n\in\R^{r'}\) and \(u_n\in\R^{N_\theta}\) such that
\begin{equation}
\label{eq:dropout-mean-decomposition}
    m_n=m_0+QD_0^{1/2}P_1b_n+u_n,
    \qquad
    \Sigma_\nu^{-1/2}F u_n\in\mathcal S^\perp.
\end{equation}
The reduced sequence \(b_n\) is exactly the same as in the no-dropout
iteration:
\begin{equation}
\label{eq:dropout-b-exact}
    b_n-b_\star=-(1-\Dt)^n\Lambda_nb_\star.
\end{equation}
In particular, the same geometric convergence bound holds:
\begin{equation}
\label{eq:dropout-b-rate}
    \|b_n-b_\star\|_2
    \leq
    \frac{\|b_\star\|_2}{\min\{1,\lambda_{\min}(\Gamma)\}}
    (1-\Dt)^n.
\end{equation}

\item \textbf{Complementary-space contraction.}
If \(\mathcal R\neq\{0\}\), then there exists \(0<\beta<1\), depending only on \(C_0,F,\Sigma_\nu,\Dt,\lambda\), such that the orthogonal observation-space error satisfies
\begin{equation}
\label{eq:dropout-h-rate}
\begin{aligned}
    \E\!\left[
        \left\|P_\perp\Sigma_\nu^{-1/2}F(m_n-m^*)\right\|^2
    \right]
    \leq
    (1-\beta)^n
    \left\|P_\perp\Sigma_\nu^{-1/2}F(m_0-m^*)\right\|^2.
\end{aligned}
\end{equation}
If \(\mathcal R=\{0\}\), the orthogonal error is identically zero. 
\item \textbf{Objective convergence.}
Furthermore, the objective error satisfies
\begin{equation}
\label{eq:dropout-objective-rate}
\begin{aligned}
    \E[\Phi_R(m_n)-\Phi_R(m^*)]
    &\leq
    \frac12\|\Gamma\|_2
    \left(
        \frac{\|b_\star\|_2}{\min\{1,\lambda_{\min}(\Gamma)\}}
    \right)^2
    (1-\Dt)^{2n}
    \\
    &\quad+
    \frac12(1-\beta)^n
    \left\|P_\perp\Sigma_\nu^{-1/2}F(m_0-m^*)\right\|^2.
\end{aligned}
\end{equation}
When \(\mathcal R=\{0\}\), the second term is zero.
\end{enumerate}
\end{theorem}

\begin{remark}[Observation subspaces and convergence mechanism]
In the linear setting, minimizing \(\Phi_R\) is equivalent
to eliminating the whitened observable error, since
\[
\Phi_R(m_n)-\Phi_R(m^*)
=
\frac12
\left\|\Sigma_\nu^{-1/2}F(m_n-m^*)\right\|^2.
\]
The observable range decomposes as
$\operatorname{Range}(\Sigma_\nu^{-1/2}F)
=
\mathcal S \oplus \mathcal R$,
where $\mathcal S$ contains the observation directions represented by the
initial ensemble and $\mathcal R$ contains the remaining observable
directions.  On this range, \(P_\perp\) extracts the component in
\(\mathcal R\).
The projected dropout correction preserves the EAKI update
in \(\mathcal S\) and supplies an additional correction
in \(\mathcal R\).
\end{remark}

\begin{proof}[Proof of \cref{thm:DropoutEAKI_ETKI}]
For $(i)$,
the particle update~\eqref{eq:dropout-particle-translation} gives
\begin{equation}
\label{eq:dropout-ensemble-translation}
\theta_{n+1}^j-m_{n+1}
=A_{n+1}\left(\widehat\theta_{n+1}^j-m_n\right),
\qquad j=1,\ldots,J.
\end{equation}
Hence \(Z_{n+1}=A_{n+1}\widehat Z_{n+1}\), exactly the anomaly update of
EAKI; translating the particles to the corrected mean does not alter
their centered anomalies. The same translation preserves the EAKI
covariance. In the linear case, the covariance recursion depends only on
the covariance, so all covariance identities in
\cref{thm:EAKI_ETKI_convergence} remain valid. In particular,
\(\Range(\widehat Z_{n+1})=\Range(Q)\), since \(D_n\succ0\).

For $(ii)$, the decomposition~\eqref{eq:augmented-Kalman-gain} expresses the Kalman gain as the sum of the original EAKI gain and the projected dropout gain.
The former equals \(\Dt C_{n+1}F^T\Sigma_\nu^{-1}\)
by~\eqref{eq:linear-gain-identity}. Thus, the mean update~\eqref{eq:dropout-mean-correction}
takes the form
\begin{subequations}\label{eq:dropout-joint-splitting}
\begin{align}
m_{n+1}&=m_n+\Dt C_{n+1}F^T\Sigma_\nu^{-1}(x-Fm_n)+d_{n+1}, \label{eq:dropout-joint-splitting-mean} \\
d_{n+1}&:=\widetilde Z_{n+1}^{\perp}(\widetilde Y_{n+1}^{\perp})^T
\left(\widetilde Y_{n+1}^{\perp}
(\widetilde Y_{n+1}^{\perp})^T+\Sigma_{\nu,n+1}\right)^{-1}(x-Fm_n). 
\end{align}
\end{subequations}
By part $(i)$,
\(\Range(\Sigma_\nu^{-1/2}\widehat Y_{n+1})=\mathcal S\).
Weighted orthogonality and
\(F\widetilde Z_{n+1}^{\perp}=\widetilde Y_{n+1}^{\perp}\)
therefore give
\(\Range(\Sigma_\nu^{-1/2}F\widetilde Z_{n+1}^{\perp})
\subseteq\mathcal R\), so
\(\Sigma_\nu^{-1/2}Fd_{n+1}\in\mathcal R\).

The decomposition~\eqref{eq:dropout-mean-decomposition}
holds at \(n=0\) with \(b_0=0\) and \(u_0=0\).
Suppose it holds at step \(n\). The accumulated
dropout displacement $u_n$ does not affect
the original EAKI increment, since
\[
    Q^TF^T\Sigma_\nu^{-1}Fu_n
    =(\Sigma_\nu^{-1/2}FQ)^T(\Sigma_\nu^{-1/2}Fu_n)=0.
\]
Substituting \eqref{eq:dropout-mean-decomposition} at $n$ into \eqref{eq:dropout-joint-splitting-mean} and using the gain identity~\eqref{eq:linear-gain-identity} together with the definition of $b_{\star}$ in~\eqref{eq:b_star} yields
\[
    b_{n+1}=b_n+\Dt\Lambda_{n+1}\Gamma(b_\star-b_n),
    \qquad u_{n+1}=u_n+d_{n+1}.
\]
The second update preserves
\(\Sigma_\nu^{-1/2}Fu_{n+1}\in\mathcal S^\perp\), completing the
induction. The first update is exactly~\eqref{eq:linear-b-recursion}, yielding
~\eqref{eq:dropout-b-exact} and~\eqref{eq:dropout-b-rate} for every sequence of masks.

For $(iii)$, 
let
\(\mathcal H_n:=\sigma(R_0,\ldots,R_{n-1})\)
denote the mask history.
Conditional expectations given \(\mathcal H_n\)
average over the newly sampled mask \(R_n\).
Define the orthogonal observation-space error
\begin{equation}
\label{eq:dropout-h-definition}
    h_n:=P_\perp\Sigma_\nu^{-1/2}F(m_n-m^*)\in\mathcal R.
\end{equation}
The inclusion \(h_n\in\mathcal R\) follows from
\(\mathcal S\subseteq\Range(\Sigma_\nu^{-1/2}F)\):
projecting a vector in this observable range onto
\(\mathcal S^\perp\) yields a vector in \(\mathcal R\). If \(\mathcal R=\{0\}\), then
~\eqref{eq:dropout-h-definition} gives \(h_n=0\) for all \(n\), and
we take \(\beta=0\). In the following, assume
\(\mathcal R\neq\{0\}\). To derive a recursion for $h_{n+1}$, we first show that the original EAKI increment in~\eqref{eq:dropout-joint-splitting-mean} does not contribute
to \(h_{n+1}\).
Indeed, \(C_{n+1}=QD_{n+1}Q^T\) and
\(P_\perp\Sigma_\nu^{-1/2}FQ=0\), so
$P_\perp\Sigma_\nu^{-1/2}F
\left[
\Dt C_{n+1}F^T\Sigma_\nu^{-1}(x-Fm_n)
\right]
=0$. 
Next, to express the contribution of dropout increment $d_{n+1}$ in~\eqref{eq:dropout-joint-splitting-mean} to $h_{n+1}$, we introduce whitened anomalies, and the Gram matrix of whitened anomalies:
\begin{equation}
\label{eq:dropout-EG-definitions}
\begin{aligned}
E_{n+1}:=\Sigma_\nu^{-1/2}\widetilde Y_{n+1}^{\perp},\quad
G_{n+1}:=E_{n+1}E_{n+1}^T.
\end{aligned}
\end{equation}
Part $(ii)$ gives \(\Range(E_{n+1})\subseteq\mathcal R\).
The first-order optimality condition
\(F^T\Sigma_\nu^{-1}(x-Fm^*)=0\) and
\(\Range(G_{n+1})=\Range(E_{n+1})\subseteq\mathcal R\) imply
\(    G_{n+1}\Sigma_\nu^{-1/2}(x-Fm_n)=-G_{n+1}h_n.
\)
Since \(\widetilde Y_{n+1}^{\perp}=\Sigma_\nu^{1/2}E_{n+1}\)
and \(\Sigma_{\nu,n+1}=\frac{1}{\Dt}\Sigma_\nu\), the covariance in the
dropout gain factors as
\[
\widetilde Y_{n+1}^{\perp}(\widetilde Y_{n+1}^{\perp})^T
+\Sigma_{\nu,n+1}
=\Sigma_\nu^{1/2}(G_{n+1}+\Dt^{-1} I)\Sigma_\nu^{1/2}.
\]
Inverting this factorization and using
\(P_\perp\Sigma_\nu^{-1/2}F\widetilde Z_{n+1}^{\perp}
=P_\perp E_{n+1}=E_{n+1}\), we can left-multiply~\eqref{eq:dropout-joint-splitting-mean} by $P_\perp\Sigma_\nu^{-1/2}F$ and then obtain
\begin{equation}
\label{eq:dropout-h-recursion}
\begin{aligned}
    h_{n+1}
    &=h_n+P_\perp\Sigma_\nu^{-1/2}Fd_{n+1}\\
    &=h_n+E_{n+1}(\widetilde Y_{n+1}^{\perp})^T \Sigma_\nu^{-1/2} (G_{n+1}+\Dt^{-1}I)^{-1}
      \Sigma_\nu^{-1/2}(x-Fm_n)\\
    &=h_n+G_{n+1}(G_{n+1}+\Dt^{-1}I)^{-1}
      \Sigma_\nu^{-1/2}(x-Fm_n)\\
    &=h_n-G_{n+1}(G_{n+1}+\Dt^{-1}I)^{-1}h_n\\
    &=(\Dt G_{n+1}+I)^{-1}h_n.
\end{aligned}
\end{equation}

By~\eqref{eq:dropout-h-recursion}, the change in squared error is
\begin{equation}
\label{eq:dropout-h-decrease}
\|h_{n+1}\|^2 -\|h_n\|^2 = -h_n^T[I-(I+\Dt G_{n+1})^{-2}]h_n. 
\end{equation} 
Next we will estimate the operator $I-(I+\Dt G_{n+1})^{-2}$ on $\mathcal{R}$ by comparing \(G_{n+1}\) with a bounded reference matrix. Since the projection $P_\perp$ only removes the original observation component,~\eqref{eq:dropout-projected-residuals} yields
\begin{equation}
\label{eq:dropout-projection-identity-E} 
E_{n+1}
=P_\perp\Sigma_\nu^{-1/2}F
  \left(R_n\widehat Z_{n+1}-\widehat Z_{n+1}L_{n+1}\right)=P_\perp\Sigma_\nu^{-1/2}FR_n\widehat Z_{n+1}, 
\end{equation}
and thus we can rewrite $G_{n+1}=E_{n+1}E_{n+1}^T$ as
\begin{equation}
G_{n+1}
=P_\perp\Sigma_\nu^{-1/2}FR_n
  (\widehat Z_{n+1}\widehat Z_{n+1}^T)
  R_nF^T\Sigma_\nu^{-1/2}P_\perp.\label{eq:dropout-projection-identity-G}
\end{equation}
Introduce the reference matrix
\begin{equation}
\label{eq:dropout-reference-Gram}
\overline G_{n+1}:=
P_\perp\Sigma_\nu^{-1/2}FR_nQQ^TR_n
F^T\Sigma_\nu^{-1/2}P_\perp.
\end{equation}
To compare it with \(G_{n+1}\), set
$\gamma^2:=
\frac{1}{\max\{\|D_0^{-1}\|_2,\|F^T\Sigma_\nu^{-1}F\|_2\}}$.
The convex-combination form of \eqref{eq:linear-precision-closed}
and \(Q^TQ=I_r\) imply \(D_n^{-1}\preceq\gamma^{-2}I_r\).
Hence the inflated prediction satisfies
\begin{equation}
\label{eq:dropout-predicted-covariance-lower}
    \widehat Z_{n+1}\widehat Z_{n+1}^T
    =Q\frac{D_n}{1-\Dt}Q^T \succeq \frac{\gamma^2}{1 - \Delta \tau}Q Q^T \succeq  \gamma^2 Q Q^T.
\end{equation}
Combining this bound with \eqref{eq:dropout-projection-identity-G}
gives \(G_{n+1}\succeq\gamma^2\overline G_{n+1}\).
Moreover, since \(\|Q\|_2=1\) and
\(\|R_n\|_2,\|P_\perp\|_2\leq1\),
\[
0\preceq\overline G_{n+1}\preceq M^2I,
\qquad M:=\|\Sigma_\nu^{-1/2}F\|_2.
\]
Spectral calculus and order reversal under matrix inversion then give
\begin{equation}
\label{eq:dropout-resolvent-comparison}
\begin{aligned}
I-(I+\Dt G_{n+1})^{-2}
\succeq
\Dt G_{n+1}(I+\Dt G_{n+1})^{-1}
\succeq
\frac{\gamma^2\Dt}{1+\gamma^2M^2\Dt}
\overline G_{n+1}.
\end{aligned}
\end{equation}

It remains to show that the reference matrix $\overline G_{n+1}$ is uniformly lower-bounded over all directions in $\mathcal{R}$, in expectation. Fix \(h\in\mathcal R\) and let
\(v:=F^T\Sigma_\nu^{-1/2}h\).
Since \(P_\perp h=h\),~\eqref{eq:dropout-reference-Gram} gives
\(h^T\overline G_{n+1}h=v^TR_nQQ^TR_nv\).
Moreover, independence of the Bernoulli coordinates gives
\(\E[R_nQQ^TR_n]=\left(\lambda^2QQ^T+
\lambda(1-\lambda)\operatorname{diag}(QQ^T)\right)\left(1-(1-\lambda)^{N_\theta}\right)^{-1}\), where the divisor $1-(1-\lambda)^{N_\theta}$ arises from the condition that $R_n\neq 0$. Combining with the coordinate-coverage condition~\eqref{eq:dropout-coordinate-coverage},
\begin{equation}  
\label{eq:dropout-estimate-hGh}
\begin{aligned}
    h^T\E[\overline G_{n+1}\mid\mathcal H_n]h
&=\E[v^TR_nQQ^TR_nv\mid\mathcal H_n]\\
&=\frac{\lambda^2\|Q^Tv\|^2+
\lambda(1-\lambda)\sum_{i=1}^{N_\theta}
(e_i^TQQ^Te_i)v_i^2}{1-(1-\lambda)^{N_\theta}}\\
&\geq\lambda(1-\lambda)\ell_Q\|v\|^2.
\end{aligned}
\end{equation}
To bound \(\|v\|\) from below by \(\|h\|\), let \(\sigma_F>0\)
be the smallest positive singular value of
\(\Sigma_\nu^{-1/2}F\). The matrix
\(\Sigma_\nu^{-1/2}FF^T\Sigma_\nu^{-1/2}\) is positive definite
on \(\Range(\Sigma_\nu^{-1/2}F)\), with smallest eigenvalue
\(\sigma_F^2\) on that range. Hence
\begin{equation}
\label{eq:dropout-automatic-observation-coercivity}
z^T\Sigma_\nu^{-1/2}FF^T\Sigma_\nu^{-1/2}z
    \geq\sigma_F^2\|z\|^2,
    \qquad \forall z\in\Range(\Sigma_\nu^{-1/2}F).
\end{equation}
Applying this to \(h\in\mathcal R\subset \Range(\Sigma_\nu^{-1/2}F) \) gives
\(\|v\|^2\geq\sigma_F^2\|h\|^2\).
Since the preceding quadratic-form estimate~\eqref{eq:dropout-estimate-hGh} holds for every
\(h\in\mathcal R\), we obtain 
\begin{equation}
\label{eq:dropout-reference-conditional-lower}
   h^T \E[\overline G_{n+1}\mid\mathcal H_n] h
    \geq\lambda(1-\lambda)\ell_Q\sigma_F^2 \|h\|^2 \qquad \forall h\in \mathcal{R}.
\end{equation}

By~\eqref{eq:dropout-h-decrease}, the squared-error decrease is
\(h_n^T[I-(I+\Dt G_{n+1})^{-2}]h_n\).
Since \(h_n\in \mathcal{R}\) is \(\mathcal H_n\)-measurable, we can apply the
preceding matrix bound~\eqref{eq:dropout-resolvent-comparison}
 inside the conditional expectation and then
use~\eqref{eq:dropout-reference-conditional-lower}:
\begin{equation}
\label{eq:dropout-resolvent-contraction}
\begin{aligned}
\E[\|h_{n+1}\|^2\mid\mathcal H_n] 
&=\|h_n\|^2-h_n^T\E[I-(I+\Dt G_{n+1})^{-2}
                 \mid\mathcal H_n]h_n\\
&\leq \|h_n\|^2-\frac{\gamma^2\Dt}{1+\gamma^2M^2\Dt}
   h_n^T\E[\overline G_{n+1}\mid\mathcal H_n]h_n\\
&\leq(1-\beta)\|h_n\|^2,
\end{aligned}
\end{equation}
where
\begin{equation}
\label{eq:dropout-beta}
\beta:=\lambda(1-\lambda)\ell_Q\sigma_F^2
       \frac{\gamma^2\Dt}{1+\gamma^2M^2\Dt}.
\end{equation}
Taking expectations in~\eqref{eq:dropout-resolvent-contraction}
and using the tower property gives
\(\E\|h_{n+1}\|^2\leq(1-\beta)\E\|h_n\|^2\).
Iterating yields
$\E\|h_n\|^2\leq(1-\beta)^n\|h_0\|^2$,
which is~\eqref{eq:dropout-h-rate} by the definition of \(h_n\)
in~\eqref{eq:dropout-h-definition}.

For $(iv)$,
the first-order optimality condition at \(m^*\) gives
\(\Phi_R(m)-\Phi_R(m^*)=\frac12\|\Sigma_\nu^{-1/2}F(m-m^*)\|^2\).
Combining this condition with the optimality condition for
the restricted minimizer
\(m_\star=m_0+QD_0^{1/2}P_1b_\star\) yields
\(\Sigma_\nu^{-1/2}F(m_\star-m^*)\perp\mathcal S\). Using~\eqref{eq:dropout-mean-decomposition} and the definition
of \(m_\star\), we write
\begin{equation}
\label{eq:mean}
m_n-m^*
    =
    QD_0^{1/2}P_1(b_n-b_\star)
    +u_n+(m_\star-m^*).
\end{equation}
Since $\Sigma_\nu^{-1/2}FQD_0^{1/2}P_1(b_n-b_\star)
    \in\mathcal S$ and
    $\Sigma_\nu^{-1/2}Fu_n,\;
    \Sigma_\nu^{-1/2}F(m_\star-m^*)
    \in\mathcal S^\perp$, applying \(P_\perp\Sigma_\nu^{-1/2}F\) to \eqref{eq:mean} and using
the definition of \(h_n\) in \eqref{eq:dropout-h-definition} gives
$h_n
    =
    \Sigma_\nu^{-1/2}Fu_n
    +\Sigma_\nu^{-1/2}F(m_\star-m^*)$.
Multiplying \eqref{eq:mean} by \(\Sigma_\nu^{-1/2}F\)
and replacing the last two terms by \(h_n\) yields
\[
    \Sigma_\nu^{-1/2}F(m_n-m^*)
    =
    \Sigma_\nu^{-1/2}FQD_0^{1/2}P_1(b_n-b_\star)+h_n.
\]
The two terms lie in \(\mathcal S\) and \(\mathcal R\),
respectively. By orthogonality and~\eqref{eq:P_1_Gamma_P1},
\begin{equation}
\label{eq:dropout-objective-identity}
\Phi_R(m_n)-\Phi_R(m^*)
=\frac12(b_n-b_\star)^T\Gamma(b_n-b_\star)
+\frac12\|h_n\|^2.
\end{equation}
Taking expectations and applying the bounds for
\(b_n-b_\star\) and \(h_n\) in
~\eqref{eq:dropout-b-rate} and~\eqref{eq:dropout-h-rate}, respectively,
proves~\eqref{eq:dropout-objective-rate}.
\end{proof}

\section{Numerical Study}
\label{sec:numerics}

In this section, we evaluate the proposed inflated projected-dropout ensemble adjustment Kalman inversion (IPD-EAKI) method on three
classes of inverse problems\footnote{Our code is available online:
\url{https://github.com/PKU-CMEGroup/InverseProblems.jl/tree/master}.}.
We compare IPD-EAKI with inflated EAKI (I-EAKI), described in
\cref{sec:EKI}, dropout ensemble Kalman inversion
(DEKI)~\cite{liu2025dropout}, and the covariance matrix adaptation
evolution strategy (CMA-ES)~\cite{hansen2003reducing}.
I-EAKI provides a baseline without the projected dropout correction,
while DEKI provides a comparison with an existing dropout method.
CMA-ES is a stochastic, derivative-free optimizer that adapts a
Gaussian search distribution.
\begin{enumerate}
    \item \textbf{Ill-conditioned linear inverse problems:} We vary the condition number of the forward operator and the number of observable directions missing from the initial ensemble subspace. 
    These experiments demonstrate how covariance inflation accelerates
    convergence within that subspace and how projected dropout
    recovers the missing observable components.
    
    \item \textbf{100-dimensional nonlinear problems:} We consider the monotone cubic, Rosenbrock, and Rastrigin problems, representing smooth nonlinear distortion, curved and ill-conditioned valleys, and highly oscillatory landscapes with multiple local minima, respectively. We assess both optimization performance and the evolution
    of ensemble spread.
    
    \item \textbf{Barotropic flow inverse problem:} We recover the initial vorticity field on the sphere from noisy zonal wind observations at two subsequent times. This example assesses the applicability of IPD-EAKI to a nonlinear PDE-constrained inverse problem.
\end{enumerate}
Throughout the experiments, I-EAKI and IPD-EAKI use a fixed artificial time step $\Delta\tau=0.5$. IPD-EAKI uses a dropout rate of $1-\lambda=0.5$. Ensemble sizes and other method-specific settings are specified for each example.

\subsection{Ill-conditioned linear inverse problems}
We consider the noiseless linear inverse problem $x=F\theta^\dagger$ with $N_\theta=N_x=100$, identity weighting $\Sigma_\nu=I$, and no prior regularization. For each condition number $\kappa\in\{1,10,10^2,10^3,10^4\}$, we construct the forward matrix as
\begin{equation}
    F=U_F\operatorname{diag}(\sigma_1,\dots,\sigma_{100})V_F^T,
    \qquad
    \sigma_i=\kappa^{-(i-1)/99},
\end{equation}
where $U_F$ and $V_F$ are random orthogonal matrices. The prescribed singular values ensure that
\(\operatorname{cond}_2(F)=\kappa\).  
We initialize $J=21$ particles in the subspace spanned by the first $20$ right singular vectors of  $F$.  Let \(Q_{\rm sub}=[v_1,\dots,v_{20}]\), where \(v_i\)
denotes the \(i\)th column of \(V_F\).  
We draw independent coefficients
\(c_0^{(j)}\sim\mathcal N(0,I_{20})\) and set
$\theta_0^{(j)}=Q_{\rm sub}c_0^{(j)}$. The initial ensemble therefore has zero empirical mean,
and its empirical covariance has rank \(20\) almost surely.
All directions in this initial subspace are observable,
since \(FQ_{\rm sub}\) has full column rank.
We construct the truth as
\begin{equation}
    \theta^\dagger
    =30\sum_{i=1}^{20+N_{\rm missing}}\frac{v_i}{\sigma_i},
\end{equation}
where \(N_{\rm missing}\in\{0,2,10\}\) is the number
of active singular directions outside the initial
ensemble subspace.
This scaling gives an observation-space amplitude of \(30\)
in each active left singular direction. 

\begin{figure}[t]
    \centering
    \includegraphics[width=0.95\textwidth]{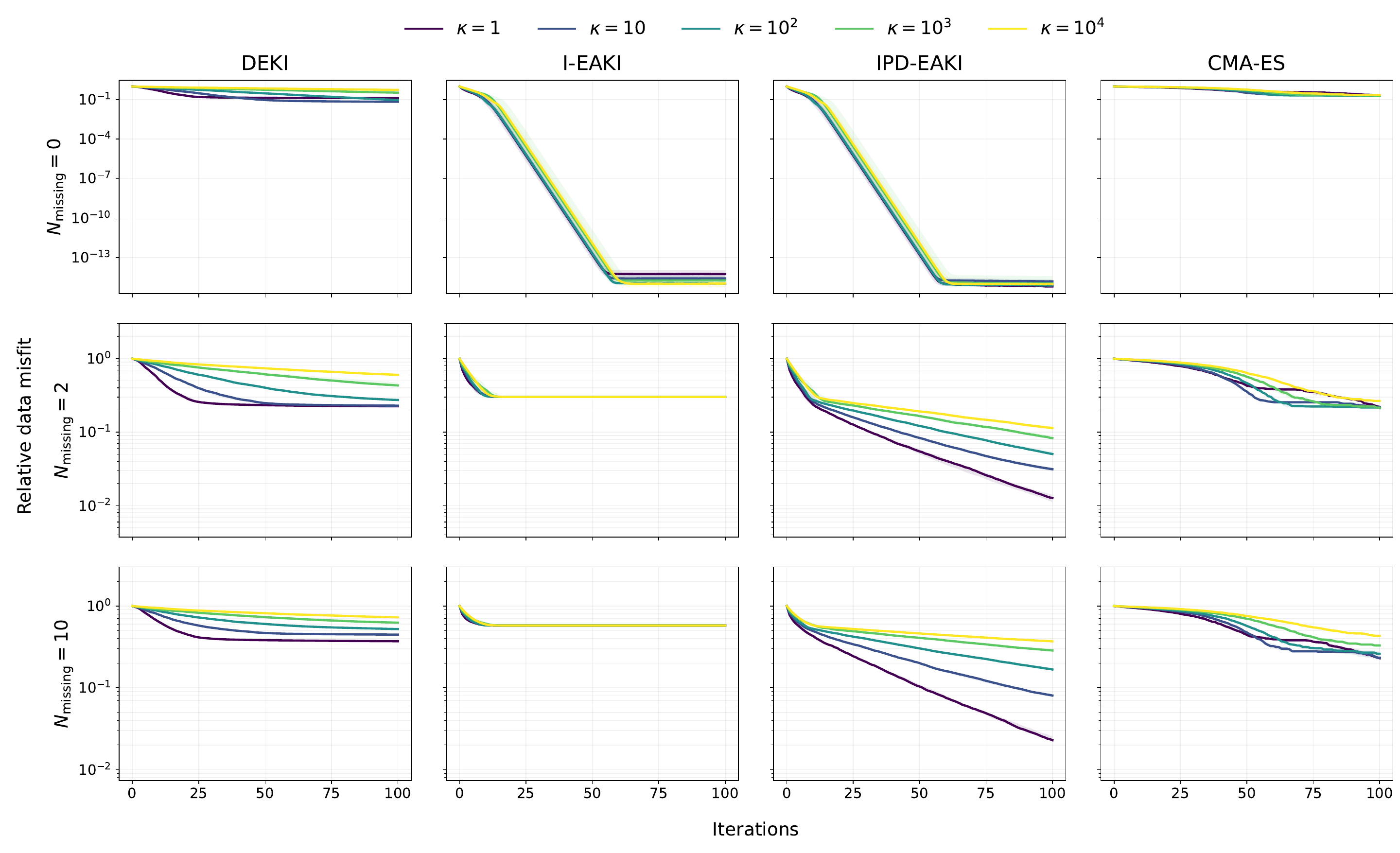}
    \caption{Condition-number sensitivity for the noiseless linear inverse problem. Columns show DEKI, I-EAKI, IPD-EAKI, and CMA-ES; rows correspond to $N_{\rm missing}=0,2,10$. Colors indicate $\kappa\in\{1,10,10^2,10^3,10^4\}$. Each curve shows the pointwise median relative data misfit over ten random seeds, with shaded regions indicating the interquartile range. The diagnostic uses the ensemble mean for the ensemble Kalman methods and the best point found so far for CMA-ES. All vertical axes are logarithmic, with different limits across rows.}
    \label{fig:cond_convergence}
\end{figure}

We compare DEKI, I-EAKI, IPD-EAKI, and CMA-ES over $100$ iterations and ten random seeds.  DEKI uses the same dropout rate and a step-size prefactor of $0.5$ for both the mean and ensemble-deviation updates. CMA-ES starts at the initial ensemble mean with step size $\sigma_0=1$ and population size $21$. Convergence is measured by the relative data misfit $\operatorname{misfit}_n=\|F m_n-x\|_2/\|x\|_2$, where $m_n$ is the ensemble mean for the ensemble Kalman methods and the best point found so far for CMA-ES.  DEKI and IPD-EAKI require approximately twice as many forward evaluations per iteration as I-EAKI and CMA-ES.

Figure~\ref{fig:cond_convergence} shows the relative data misfit
as a function of iteration.
For \(N_{\rm missing}=0\), I-EAKI and IPD-EAKI exhibit nearly
identical decay to numerical precision, consistent with the
asymptotic factor \(1-\Delta\tau\), which is independent of
the condition number \(\kappa\); see
\cref{thm:EAKI_ETKI_convergence}.
For \(N_{\rm missing}>0\), I-EAKI approaches the subspace
misfit floor, whereas the projected dropout mean correction
enables IPD-EAKI to reduce the misfit further, although this
additional reduction slows as \(\kappa\) increases.
Under the same iteration budget, IPD-EAKI achieves lower
final misfits than DEKI and CMA-ES, both of which remain
sensitive to conditioning.

\subsection{100-dimensional nonlinear problems}
\label{sec:numerical_nonlinear}
We consider three nonlinear least-squares problems in
\(N_\theta=100\) dimensions:
$\min 
\Phi(\theta)$,
$\Phi(\theta) = \frac12\|\mathcal F(\theta)\|_2^2$.
Each problem has a known global minimizer \(\theta^\ast\)
satisfying \(\mathcal F(\theta^\ast)=0\).

\textbf{(1) Monotone cubic problem.}
The forward map is
\begin{equation}
\mathcal F(\theta)
=
g(B\theta)-g(B\theta^\ast) \in \R^{N_\theta},  
\end{equation}
where $g(z)=z+ 2 z^3$ is applied componentwise, $B \in \R^{N_\theta \times N_\theta}$ is a random orthogonal matrix, 
and $\theta^\ast$ is randomly generated.
Since $g'(z)>0$, this example has a unique global minimizer
and tests smooth nonlinear distortion without spurious local minima.

\textbf{(2) Rosenbrock problem.}
The classical Rosenbrock objective admits the following
residual representation:
\[
\mathcal F(\theta)=[f_1(\theta),\cdots,f_{2N_\theta-2}(\theta)]^T 
\quad 
f_{2i-1}(\theta)=\sqrt{200}
(\theta_{i+1}-\theta_i^2),
\quad
f_{2i}(\theta)=\sqrt{2}(1-\theta_i),
\]
The unique global minimizer is $\theta^\ast=(1,\ldots,1)^T$.
This problem tests convergence in curved,
ill-conditioned valleys.

\textbf{(3) Rastrigin problem.}
The classical Rastrigin objective admits the residual
representation
\[
\mathcal F(\theta)
=
\sqrt{2}
\begin{bmatrix}
\theta\\
\sqrt{20}\sin(\pi\theta)
\end{bmatrix} \in \R^{2N_\theta},
\]
where the sine function is applied componentwise.
The unique global minimizer is $\theta^\ast=0$.
This problem tests robustness on a highly oscillatory
landscape with many local minima.

We compare DEKI, IPD-EAKI, and CMA-ES using the same budget of forward evaluations and ten random seeds. 
For all experiments, we use \(J\in\{20,50\}\) particles. Following the standard setup,
we sample the initial ensembles from Gaussian distributions
\(\N(0,I)\) for the monotone cubic and Rosenbrock
problems, and uniformly from
\([-3,5]^{N_\theta}\) for the Rastrigin problem.

\begin{figure}[htbp]
    \centering
    \includegraphics[width=0.95\textwidth]{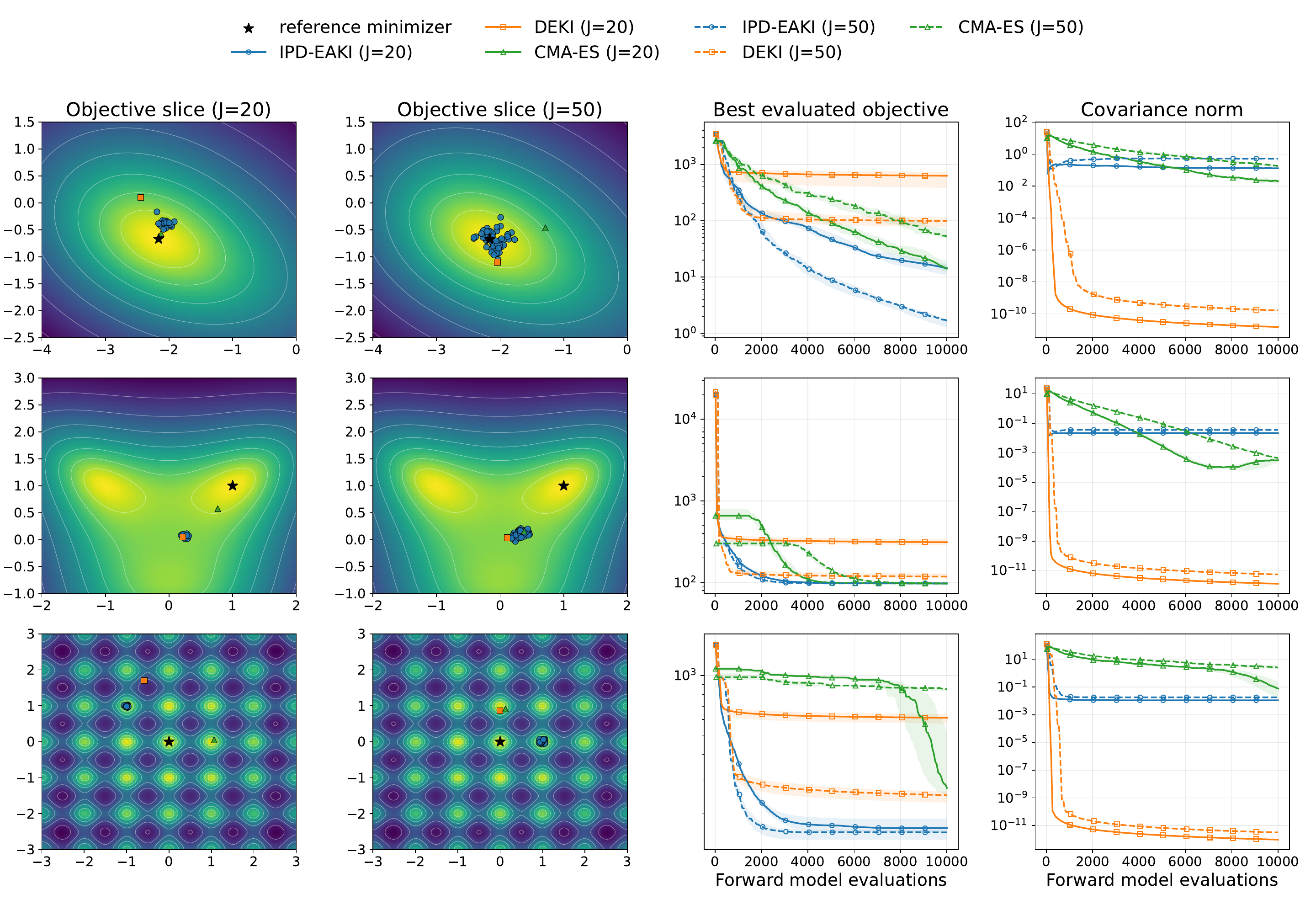}
    \caption{
    Nonlinear results on the monotone cubic, Rosenbrock, and Rastrigin problems, shown from top to bottom. The first two columns show two-dimensional objective slices with the final particle positions for \(J=20\) and \(J=50\), respectively, from one representative run. The third and fourth columns show the best evaluated objective value and covariance norm against forward model evaluations. Colors distinguish IPD-EAKI, DEKI, and CMA-ES; solid and dashed curves correspond to \(J=20\) and \(J=50\). Curves show the median over 10 independent random seeds, with shaded regions indicating the interquartile range.
    }
    \label{fig:nonlinear_results}
\end{figure}

\cref{fig:nonlinear_results} shows representative two-dimensional objective slices with projected final particle locations, and the best evaluated objective and ensemble covariance norm. The slices are obtained by varying the first two parameter coordinates $\theta_1,\theta_2$ while fixing all remaining coordinates at their reference minimizer values $\theta^*_3,\theta^*_4,\cdots,\theta^*_{N_\theta}$. Across all three problems, IPD-EAKI and DEKI initially reduce the objective faster than CMA-ES, but DEKI subsequently stagnates, coinciding with
rapid covariance collapse. IPD-EAKI achieves lower
median best objective values while maintaining
ensemble diversity.
On Rosenbrock, IPD-EAKI rapidly reaches objective values comparable to CMA-ES, although both methods subsequently progress
slowly, consistent with the difficulty of navigating the
narrow, curved valley. 
On Rastrigin, IPD-EAKI also outperforms CMA-ES, suggesting
that sustained ensemble spread may facilitate exploration
of the multimodal landscape. 

\subsection{Barotropic flow inverse problem}
\label{sec:dropout-barotropic}

\begin{figure}[htbp]
    \centering
    \includegraphics[width=0.9\linewidth]{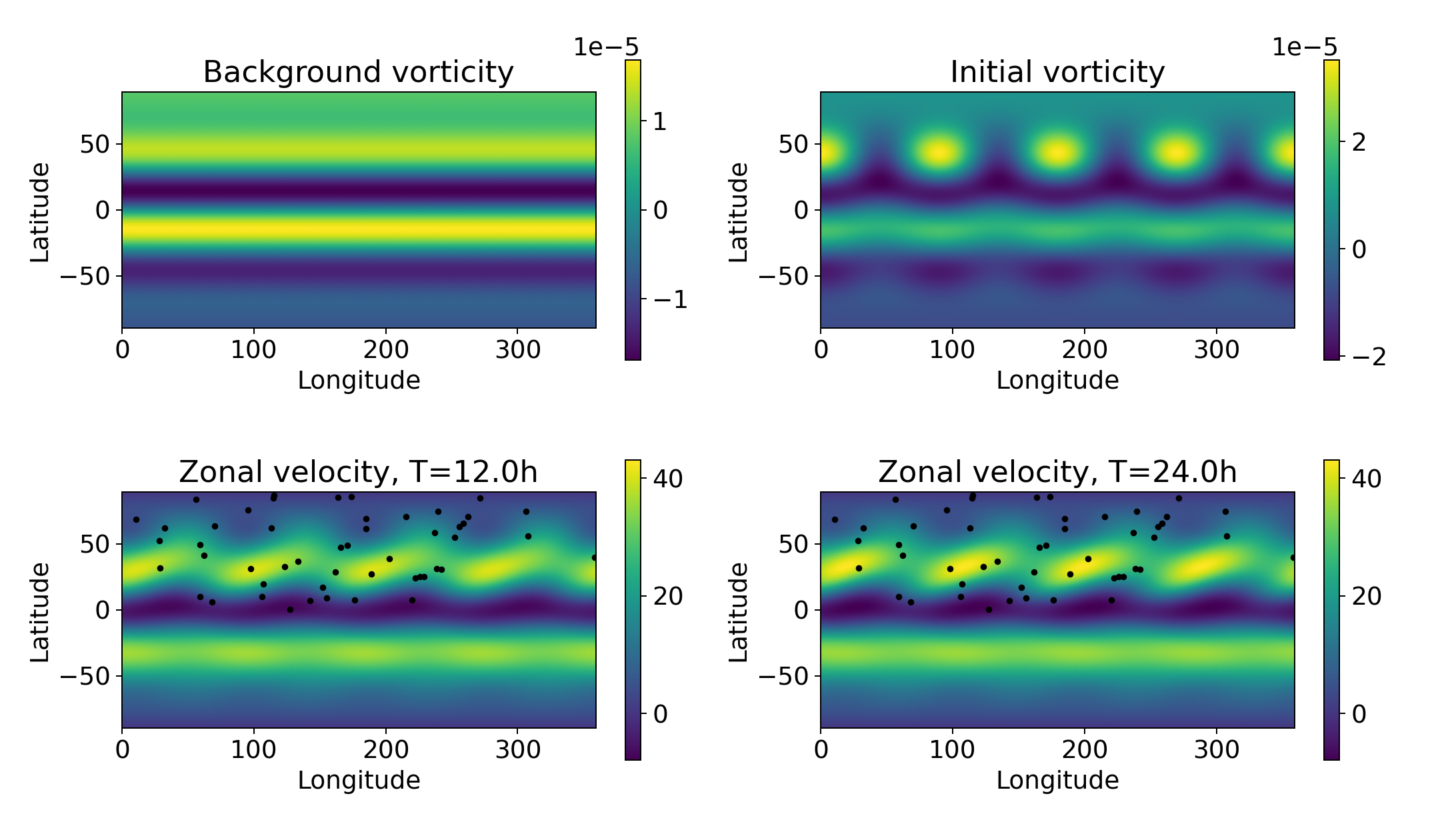}
    \caption{The upper panels show the background vorticity (left) and
    the full initial vorticity (right). The lower panels show the zonal
    wind at $T=12.0\,\mathrm{h}$ (left) and $T=24.0\,\mathrm{h}$ (right).
    Black dots mark the same 50 observation locations in the Northern Hemisphere.}
    \label{fig:vor_0}
\end{figure}

Finally, we consider the incompressible Euler equations on a
rotating sphere, written in vorticity form as
\begin{equation}
\label{eq:barotropic}
    \partial_t\zeta=-v\cdot\nabla(\zeta+f),
    \qquad \nabla^2\psi=\zeta,
    \qquad v=k\times\nabla\psi,
\end{equation}
where $\zeta$ is relative vorticity, $\psi$ is the streamfunction, $k$ is the
radial unit vector, and $f=2\Omega\sin\varphi$ is the Coriolis parameter at latitude $\varphi$.
We use Earth's radius $R=6.3712\times10^6\,\mathrm{m}$ and rotation rate
$\Omega=7.292\times10^{-5}\,\mathrm{s}^{-1}$.
Following an idealized configuration for Northern Hemisphere
winter flow in the upper troposphere~\cite{held1985pseudomomentum,held1987linear},
we construct the reference initial vorticity as \(\zeta_{\rm ref}=\zeta_b+\zeta'\).
Here, \(\zeta_b\) is the background vorticity associated with
the zonal wind
\begin{equation}
\label{eq:dropout-u-b}
    u_b(\varphi)=25\cos\varphi-30\cos^3\varphi
        +300\sin^2\varphi\cos^6\varphi
        \quad \mathrm{m}\,\mathrm{s}^{-1},
\end{equation}
and the disturbance is given by
$\zeta'(\ell,\varphi)=\frac{A}{2}\cos\varphi
    \exp\!\left[-\left(\frac{\varphi-\varphi_0}{\varphi_W}\right)^2\right]
    \cos(4\ell)$,
where $\ell$ is longitude, $\varphi_0=45^\circ\mathrm{N}$,
$\varphi_W=15^\circ$, and $A=8.0\times10^{-5}\,\mathrm{s}^{-1}$.
The initial condition is
shown in \cref{fig:vor_0}.

The forward problem is solved at T85 spectral resolution on a
$512\times256$ longitude--latitude grid, using a filtered leapfrog
integrator \cite{held1994proposal} with a $1800\,\mathrm{s}$ time step.

For the inverse problem, we recover the initial vorticity field $\zeta_0$ from
noisy pointwise observations of the zonal wind velocity.
Observations are collected at the same 50 randomly selected
locations in the Northern Hemisphere at \(T=12\,\mathrm{h}\)
and \(T=24\,\mathrm{h}\) and perturbed by Gaussian noise with
covariance matrix \(\Sigma_\eta\) determined by componentwise standard
deviations equal to \(5\%\) of the absolute values of the
corresponding reference data, yielding the observation vector
\(y_{\rm obs}\) with \(N_y=100\) entries
(see the bottom row of \cref{fig:vor_0}).
The background vorticity
$\zeta_b$ is used as the prior mean of the full initial field.
The zero-mean real vorticity field on the sphere can be approximated with $N_{\text{trunc}} = 7$ truncated wavenumber with triangular truncation, as follows,
\begin{equation}
\label{eq:sph-harmonics}
    \zeta_0(\ell, \varphi) = \zeta_b + \sum_{n=1}^{N_{\text{trunc}}} a_{0,n} P_{0,n}(\varphi) + \sum_{m = 1}^{N_{\text{trunc}}} \sum_{n=m}^{N_{\text{trunc}}} 2P_{m,n}(\varphi) \left(a^c_{m,n} \cos(m\ell) +  a^s_{m,n} \sin(m\ell) \right). \\
\end{equation}
Collecting these coefficients gives $a\in\mathbb{R}^{63}$. 
Solving the~\eqref{eq:barotropic} gives the forward model mapping $\mathcal{G}: a \rightarrow y_{\rm obs}$. The Gaussian
prior $\mathcal{N}(0, \Sigma_0)$ with $\Sigma_0=9I_{63}$ is incorporated through the augmented map
\begin{equation}
\label{eq:barotropic-augmented-map}
    {\mathcal{F}}(a)=
    \begin{bmatrix}
        \mathcal{G}(a)\\
        a
    \end{bmatrix},
    \qquad
    x=
    \begin{bmatrix}
        y_{\mathrm{obs}}\\
        0
    \end{bmatrix},
    \qquad 
    \Sigma_{\nu} = \begin{bmatrix}
        \Sigma_{\rm \eta} & 0 \\
        0 & \Sigma_0
    \end{bmatrix}.
\end{equation}

\begin{figure}[htbp]
    \centering
    \includegraphics[width=0.8\linewidth]{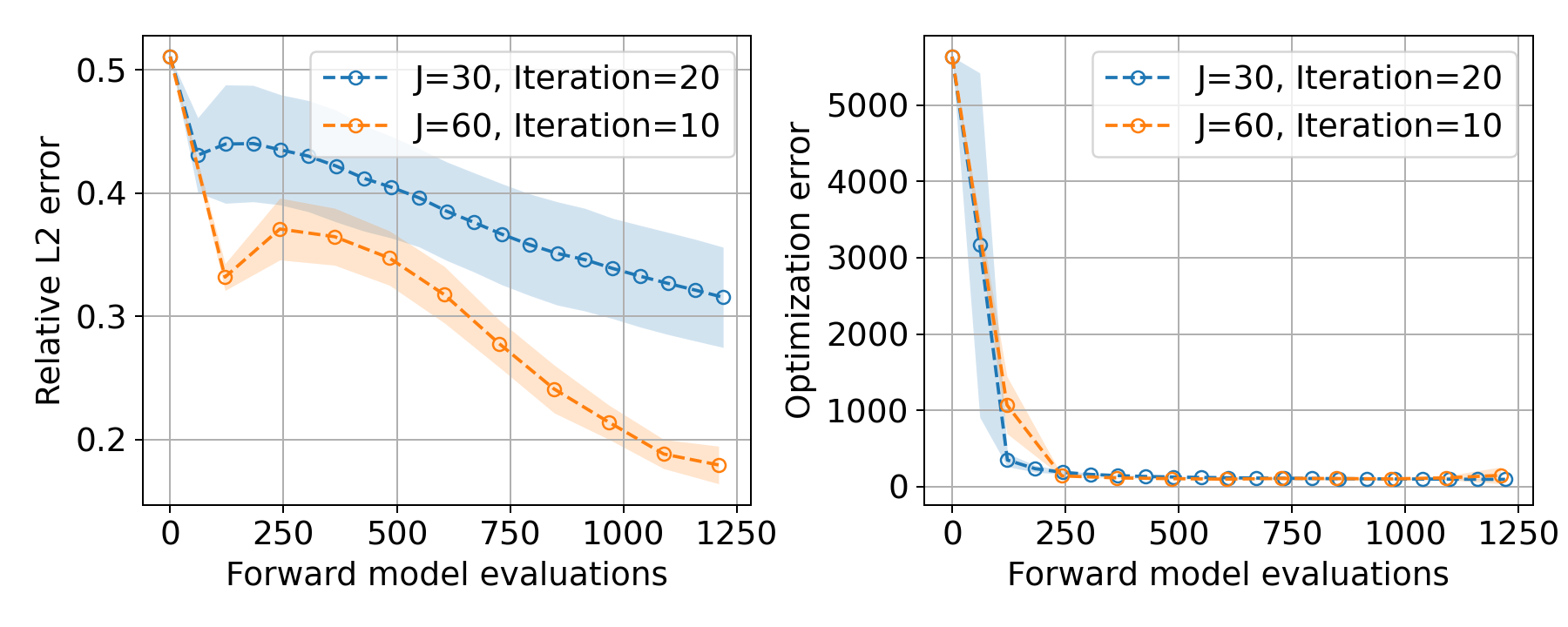}
    \caption{Relative vorticity error $\frac{\left\|\zeta_0(m_n)-\zeta_{\mathrm{ref}}\right\|_2}
    {\left\|\zeta_{\mathrm{ref}}\right\|_2}$ (left) and regularized objective $\Phi_R(m_n)$ (right) for the two IPD-EAKI configurations. Each curve shows the mean over 10 runs with different dropout seeds, keeping the problem, observations, and initialization fixed within each configuration. Shaded regions indicate the mean $\pm$ sample standard deviation.}
    \label{fig:barotropic-results}
\end{figure}

\begin{figure}[htbp]
    \centering
    \includegraphics[width=0.98\linewidth,trim={0mm 10mm 0mm 8mm},clip]{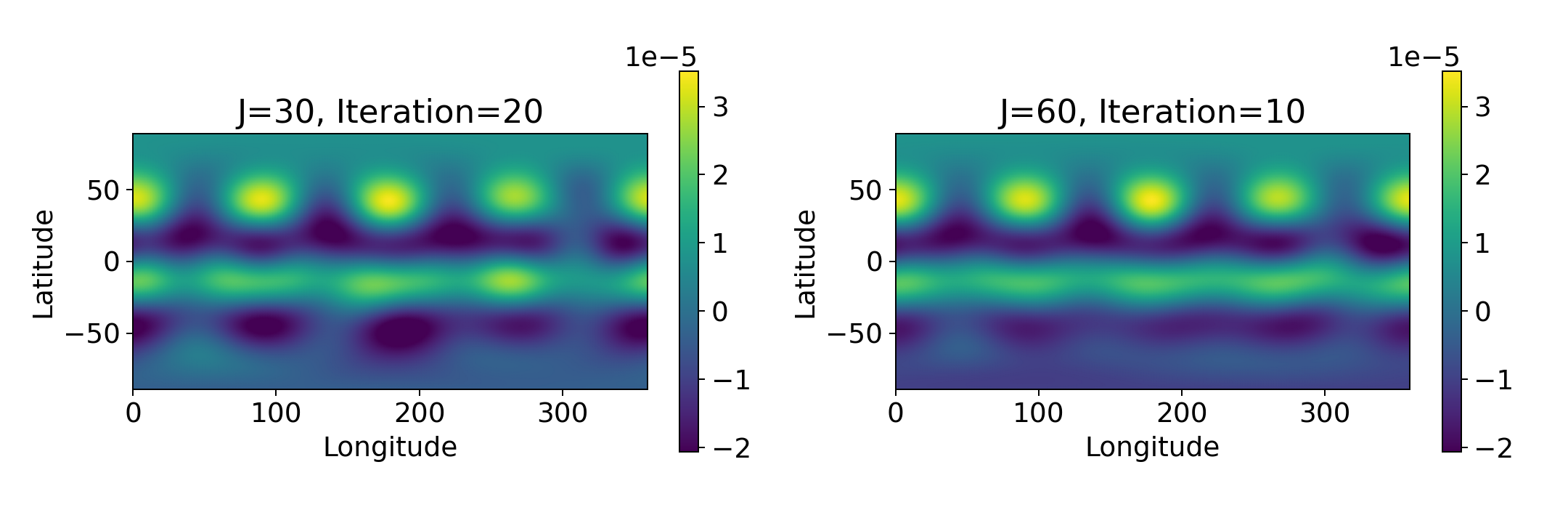}\\
    \includegraphics[width=0.98\linewidth,trim={0mm 10mm 0mm 8mm},clip]{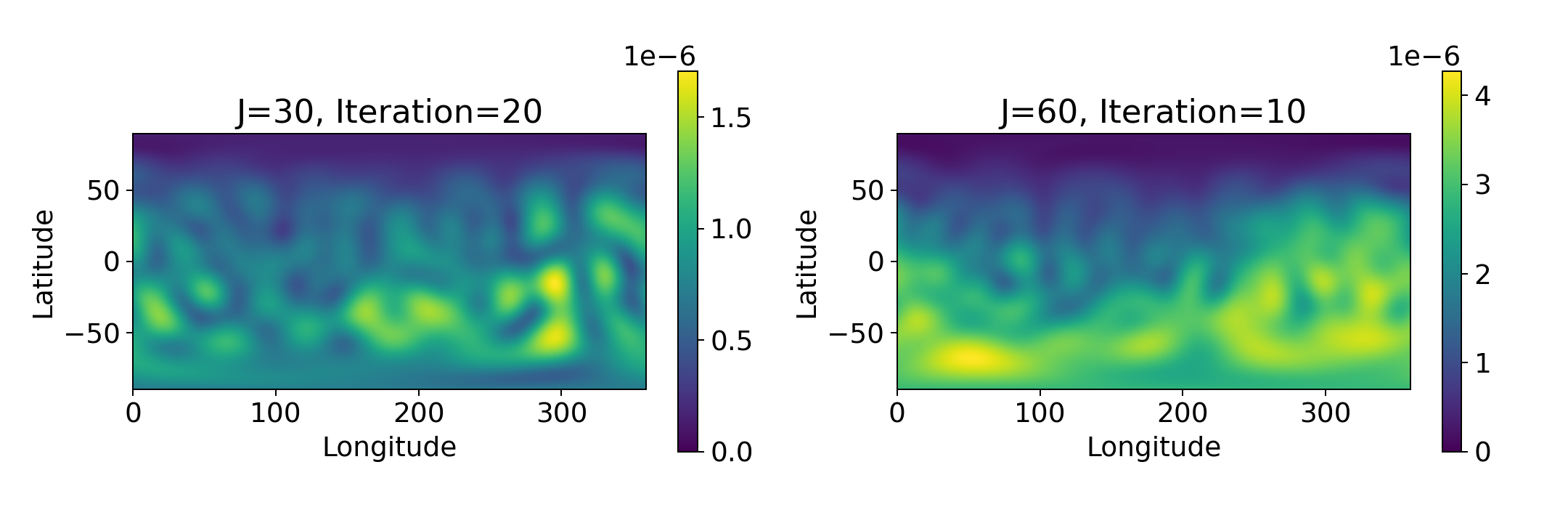}
    \caption{Reconstructed initial vorticity fields from the final ensemble means (top) and pointwise standard deviations over the final ensemble members (bottom), using the first of 10 trials for each configuration: $(J,N_{\rm iter})=(30,20)$ (left) and $(60,10)$ (right). The reference initial vorticity is shown in \cref{fig:vor_0}.}
    \label{fig:barotropic-vorticity}
\end{figure}

We compare two IPD-EAKI configurations,
\((J,N_{\rm iter})=(30,20)\) and \((60,10)\),
using comparable budgets of forward evaluations.
Each configuration is run 10 times with different dropout seeds,
keeping the problem, observations, and initialization fixed.
\cref{fig:barotropic-results} shows the relative vorticity
error and the regularized objective \(\Phi_R\) associated
with the augmented formulation \eqref{eq:barotropic-augmented-map}.
Both configurations reduce these quantities, with the
60-member ensemble achieving lower values.
This improvement may reflect the larger initial ensemble
subspace, which provides more directions for direct EAKI
updates, potentially accelerating convergence and reducing
reliance on projected dropout for exploration beyond that
subspace.
The reconstructed initial vorticity fields in \cref{fig:barotropic-vorticity}-top
also show that the 60-member ensemble provides a more
accurate reconstruction.
The final ensemble standard-deviation fields in
\cref{fig:barotropic-vorticity}-bottom highlight estimated uncertainty
in the Southern Hemisphere, where no direct observations
are available.
The higher standard deviations obtained with \(J=60\) may
reflect the broader range of parameter directions represented
by its ensemble covariance.

\section{Conclusion}
\label{sec:conclusion}

We enhance the EKI framework by combining multiplicative covariance inflation with projected dropout to improve exploration and convergence while retaining its derivative-free, parallel structure. Numerical experiments on ill-conditioned linear systems, high-dimensional nonlinear benchmarks, and barotropic flow inversion support the complementary roles of the two mechanisms.
Future work will investigate adaptive parameter selection and line-search or trust-region acceptance criteria to further improve robustness in nonlinear problems.

\paragraph{Acknowledgments} 
We acknowledge funding support from National Key R\&D Program of China 2025YFA1018700, 
National Natural Science Foundation of China (No.12471403, 62595771, and 12288101), the Fundamental Research Funds for the Central Universities of China, and the support of the high-performance computing platform of Peking University.

\paragraph{AI disclosure.} 
We used GPT-6 (OpenAI) to assist in developing theoretical arguments
for the linear case. We subsequently verified, simplified, and
organized these arguments into the proofs of
\cref{thm:EAKI_ETKI_convergence,thm:DropoutEAKI_ETKI}.
GPT-6 was also used to improve the language. The authors take full responsibility for the
final content.

\section*{References}
\bibliographystyle{unsrt}
\bibliography{references}

@article{vernon2025nesterov,
  title = {{Nesterov} acceleration for ensemble {Kalman} inversion and variants},
  author = {Vernon, Sydney and Bach, Eviatar and Dunbar, Oliver R. A.},
  journal = {Journal of Computational Physics},
  volume = {535},
  pages = {114063},
  year = {2025},
  publisher = {Elsevier},
  doi = {10.1016/j.jcp.2025.114063}
}

@article{huang2022efficient,
  title={Efficient derivative-free Bayesian inference for large-scale inverse problems},
  author={Huang, Daniel Zhengyu and Huang, Jiaoyang and Reich, Sebastian and Stuart, Andrew M},
  journal={Inverse Problems},
  volume={38},
  number={12},
  pages={125006},
  year={2022},
  publisher={IOP Publishing}
}

@article{held1994proposal,
  title={A proposal for the intercomparison of the dynamical cores of atmospheric general circulation models},
  author={Held, Isaac M and Suarez, Max J},
  journal={Bulletin of the American Meteorological Society},
  volume={75},
  number={10},
  pages={1825--1830},
  year={1994},
  publisher={American Meteorological Society}
}

@article{chada2020tikhonov,
  title={Tikhonov regularization within ensemble {K}alman inversion},
  author={Chada, Neil K and Stuart, Andrew M and Tong, Xin T},
  journal={SIAM Journal on Numerical Analysis},
  volume={58},
  number={2},
  pages={1263--1294},
  year={2020},
  publisher={SIAM}
}

@article{liu2025dropout,
  title={Dropout ensemble Kalman inversion for high dimensional inverse problems},
  author={Liu, Shuigen and Reich, Sebastian and Tong, Xin T},
  journal={SIAM Journal on Numerical Analysis},
  volume={63},
  number={2},
  pages={685--715},
  year={2025},
  publisher={SIAM}
}

@article{ding2021ensemble,
  title={Ensemble Kalman inversion: mean-field limit and convergence analysis},
  author={Ding, Zhiyan and Li, Qin},
  journal={Statistics and computing},
  volume={31},
  number={1},
  pages={9},
  year={2021},
  publisher={Springer}
}

@article{gao2021bi,
  title={A bi-fidelity ensemble {K}alman method for {PDE}-constrained inverse problems in computational mechanics},
  author={Gao, Han and Wang, Jian-Xun},
  journal={Computational Mechanics},
  volume={67},
  number={4},
  pages={1115--1131},
  year={2021},
  publisher={Springer}
}

@article{stuart2010inverse,
  title={Inverse problems: {A} {B}ayesian perspective},
  author={Stuart, Andrew M},
  journal={Acta numerica},
  volume={19},
  pages={451--559},
  year={2010},
  publisher={Cambridge University Press}
}

@article{tippett2003ensemble,
  title={Ensemble square root filters},
  author={Tippett, Michael K and Anderson, Jeffrey L and Bishop, Craig H and Hamill, Thomas M and Whitaker, Jeffrey S},
  journal={Monthly Weather Review},
  volume={131},
  number={7},
  pages={1485--1490},
  year={2003}
}

@article{UKI,
  title={Iterated Kalman Methodology For Inverse Problems},
  author={Huang, Daniel Zhengyu and Schneider, Tapio and Stuart, Andrew M},
  journal={Journal of Computational Physics},
  pages={111262},
  year={2022},
  publisher={Elsevier}
}

@article{kalman1960new,
  title={A new approach to linear filtering and prediction problems},
  author={{K}alman, Rudolph Emil},
  journal={J. Basic Eng. Mar},
  volume={82},
  number={1},
  pages={35-45},
  year={1960}
}

@article{evensen1994sequential,
  title={Sequential data assimilation with a nonlinear quasi-geostrophic model using {M}onte {C}arlo methods to forecast error statistics},
  author={Evensen, Geir},
  journal={Journal of Geophysical Research: Oceans},
  volume={99},
  number={C5},
  pages={10143--10162},
  year={1994},
  publisher={Wiley Online Library}
}

@article{parzer2022convergence,
  title={On convergence rates of adaptive ensemble Kalman inversion for linear ill-posed problems},
  author={Parzer, Fabian and Scherzer, Otmar},
  journal={Numerische Mathematik},
  volume={152},
  number={2},
  pages={371--409},
  year={2022},
  publisher={Springer}
}

@article{chada2022convergence,
  title={Convergence acceleration of ensemble Kalman inversion in nonlinear settings},
  author={Chada, Neil and Tong, Xin},
  journal={Mathematics of Computation},
  volume={91},
  number={335},
  pages={1247--1280},
  year={2022}
}

@article{bungert2023complete,
  title={Complete deterministic dynamics and spectral decomposition of the linear ensemble Kalman inversion},
  author={Bungert, Leon and Wacker, Philipp},
  journal={SIAM/ASA Journal on Uncertainty Quantification},
  volume={11},
  number={1},
  pages={320--357},
  year={2023},
  publisher={SIAM}
}

@article{schillings2017analysis,
  title={Analysis of the ensemble {K}alman filter for inverse problems},
  author={Schillings, Claudia and Stuart, Andrew M},
  journal={SIAM Journal on Numerical Analysis},
  volume={55},
  number={3},
  pages={1264--1290},
  year={2017},
  publisher={SIAM}
}

@article{schillings2018convergence,
  title={Convergence analysis of ensemble {K}alman inversion: {T}he linear, noisy case},
  author={Schillings, Claudia and Stuart, Andrew M},
  journal={Applicable Analysis},
  volume={97},
  number={1},
  pages={107--123},
  year={2018},
  publisher={Taylor \& Francis}
}

@article{chen2012ensemble,
  title={Ensemble randomized maximum likelihood method as an iterative ensemble smoother},
  author={Chen, Yan and Oliver, Dean S},
  journal={Mathematical Geosciences},
  volume={44},
  number={1},
  pages={1--26},
  year={2012},
  publisher={Springer}
}

@article{emerick2013investigation,
  title={Investigation of the sampling performance of ensemble-based methods with a simple reservoir model},
  author={Emerick, Alexandre A and Reynolds, Albert C},
  journal={Computational Geosciences},
  volume={17},
  number={2},
  pages={325--350},
  year={2013},
  publisher={Springer}
}

@article{iglesias2013ensemble,
  title={Ensemble {K}alman methods for inverse problems},
  author={Iglesias, Marco A and Law, Kody JH and Stuart, Andrew M},
  journal={Inverse Problems},
  volume={29},
  number={4},
  pages={045001},
  year={2013},
  publisher={IOP Publishing}
}

@inproceedings{shin2024multi,
  title={Multi-fidelity modeling and uncertainty quantification of heterogeneous roughness},
  author={Shin, YoungIn and Chan, Miles and Wang, Jianyu and Zahtila, Tony and Gorle, Catherine and Iaccarino, Gianluca and Howland, Michael},
  booktitle={Center for Turbulence Research Proceedings of the Summer Program},
  year={2024}
}

@article{zhang2022ensemble,
  title={Ensemble Kalman method for learning turbulence models from indirect observation data},
  author={Zhang, Xin-Lei and Xiao, Heng and Luo, Xiaodong and He, Guowei},
  journal={Journal of Fluid Mechanics},
  volume={949},
  pages={A26},
  year={2022},
  publisher={Cambridge University Press}
}

@article{lea2000sensitivity,
  title={Sensitivity analysis of the climate of a chaotic system},
  author={Lea, Daniel J and Allen, Myles R and Haine, Thomas WN},
  journal={Tellus A: Dynamic Meteorology and Oceanography},
  volume={52},
  number={5},
  pages={523--532},
  year={2000},
  publisher={Taylor \& Francis}
}

@article{wang2014least,
  title={Least squares shadowing sensitivity analysis of chaotic limit cycle oscillations},
  author={Wang, Qiqi and Hu, Rui and Blonigan, Patrick},
  journal={Journal of Computational Physics},
  volume={267},
  pages={210--224},
  year={2014},
  publisher={Elsevier}
}

@article{kovachki2019ensemble,
  title={Ensemble {K}alman inversion: {A} derivative-free technique for machine learning tasks},
  author={Kovachki, Nikola B and Stuart, Andrew M},
  journal={Inverse Problems},
  volume={35},
  number={9},
  pages={095005},
  year={2019},
  publisher={IOP Publishing}
}

@article{anderson2001ensemble,
  title={An ensemble adjustment {K}alman filter for data assimilation},
  author={Anderson, Jeffrey L},
  journal={Monthly weather review},
  volume={129},
  number={12},
  pages={2884--2903},
  year={2001}
}

@article{bishop2001adaptive,
  title={Adaptive sampling with the ensemble transform {K}alman filter. {P}art {I}: {T}heoretical aspects},
  author={Bishop, Craig H and Etherton, Brian J and Majumdar, Sharanya J},
  journal={Monthly weather review},
  volume={129},
  number={3},
  pages={420--436},
  year={2001}
}

@article{chen2023sampling,
  title = {Sampling via gradient flows in the space of probability measures},
  author = {Chen, Yifan and Huang, Daniel Zhengyu and Huang, Jiaoyang and Reich, Sebastian and Stuart, Andrew M},
  journal = {Mathematics of Computation},
  year = {2026},
  doi = {10.1090/mcom/4186}
}

@article{law2015data,
  title={Data assimilation},
  author={Law, Kody and Stuart, Andrew and Zygalakis, Kostas},
  journal={Cham, Switzerland: Springer},
  year={2015},
  publisher={Springer}
}

@inproceedings{chopin2023connection,
  title = {A connection between tempering and entropic mirror descent},
  author = {Chopin, Nicolas and Crucinio, Francesca and Korba, Anna},
  booktitle = {Proceedings of the 41st International Conference on Machine Learning},
  series = {Proceedings of Machine Learning Research},
  volume = {235},
  pages = {8782--8800},
  year = {2024},
  publisher = {PMLR},
  eprint = {2310.11914},
  archiveprefix = {arXiv},
  url = {https://proceedings.mlr.press/v235/chopin24a.html}
}

@article{dunbar2022ensemble,
  title={Ensemble inference methods for models with noisy and expensive likelihoods},
  author={Dunbar, Oliver RA and Duncan, Andrew B and Stuart, Andrew M and Wolfram, Marie-Therese},
  journal={SIAM Journal on Applied Dynamical Systems},
  volume={21},
  number={2},
  pages={1539--1572},
  year={2022},
  publisher={SIAM}
}

@article{emerick2013ensemble,
  title={Ensemble smoother with multiple data assimilation},
  author={Emerick, Alexandre A and Reynolds, Albert C},
  journal={Computers \& Geosciences},
  volume={55},
  pages={3--15},
  year={2013},
  publisher={Elsevier}
}

@article{dunbar2021calibration,
  title={Calibration and uncertainty quantification of convective parameters in an idealized GCM},
  author={Dunbar, Oliver RA and Garbuno-Inigo, Alfredo and Schneider, Tapio and Stuart, Andrew M},
  journal={Journal of Advances in Modeling Earth Systems},
  volume={13},
  number={9},
  pages={e2020MS002454},
  year={2021},
  publisher={Wiley Online Library}
}

@book{kaipio2006statistical,
  title={Statistical and computational inverse problems},
  author={Kaipio, Jari and Somersalo, Erkki},
  volume={160},
  year={2006},
  publisher={Springer Science \& Business Media}
}

@incollection{dashti2013bayesian,
  title = {The {B}ayesian approach to inverse problems},
  author = {Dashti, Masoumeh and Stuart, Andrew M},
  editor = {Ghanem, Roger and Higdon, David and Owhadi, Houman},
  booktitle = {Handbook of Uncertainty Quantification},
  pages = {311--428},
  year = {2017},
  publisher = {Springer},
  address = {Cham},
  doi = {10.1007/978-3-319-12385-1_7}
}

@book{doucet2001sequential,
  editor    = {Doucet, Arnaud and de Freitas, Nando and Gordon, Neil},
  title     = {Sequential {Monte Carlo} Methods in Practice},
  publisher = {Springer},
  address   = {New York, NY},
  year      = {2001},
  doi       = {10.1007/978-1-4757-3437-9}
}

@article{reich2013nonparametric,
  title={A nonparametric ensemble transform method for {B}ayesian inference},
  author={Reich, Sebastian},
  journal={SIAM Journal on Scientific Computing},
  volume={35},
  number={4},
  pages={A2013--A2024},
  year={2013},
  publisher={SIAM}
}

@article{pidstrigach2023affine,
  title={Affine-invariant ensemble transform methods for logistic regression},
  author={Pidstrigach, Jakiw and Reich, Sebastian},
  journal={Foundations of Computational Mathematics},
  volume={23},
  number={2},
  pages={675--708},
  year={2023},
  publisher={Springer}
}

@article{ruchi2019transform,
  title={Transform-based particle filtering for elliptic Bayesian inverse problems},
  author={Ruchi, Sangeetika and Dubinkina, Svetlana and Iglesias, Marco A},
  journal={Inverse Problems},
  volume={35},
  number={11},
  pages={115005},
  year={2019},
  publisher={IOP Publishing}
}

@article{weissmann2022gradient,
  title={Gradient flow structure and convergence analysis of the ensemble Kalman inversion for nonlinear forward models},
  author={Weissmann, Simon},
  journal={Inverse Problems},
  volume={38},
  number={10},
  pages={105011},
  year={2022},
  publisher={IOP Publishing}
}

@article{baptista2024representation,
  title={On the representation and learning of monotone triangular transport maps},
  author={Baptista, Ricardo and Marzouk, Youssef and Zahm, Olivier},
  journal={Foundations of Computational Mathematics},
  volume={24},
  number={6},
  pages={2063--2108},
  year={2024},
  publisher={Springer}
}

@article{marzouk2017sampling,
  title={Sampling via measure transport: An introduction},
  author={Marzouk, Youssef and Moselhy, Tarek and Parno, Matthew and Spantini, Alessio},
  journal={Springer Books},
  pages={785--825},
  year={2017},
  publisher={Springer}
}

@article{del2006sequential,
  title={Sequential monte carlo samplers},
  author={Del Moral, Pierre and Doucet, Arnaud and Jasra, Ajay},
  journal={Journal of the Royal Statistical Society Series B: Statistical Methodology},
  volume={68},
  number={3},
  pages={411--436},
  year={2006},
  publisher={Oxford University Press}
}

@article{gu2007iterative,
  title={An iterative ensemble Kalman filter for multiphase fluid flow data assimilation},
  author={Gu, Yaqing and Oliver, Dean S},
  journal={Spe Journal},
  volume={12},
  number={04},
  pages={438--446},
  year={2007},
  publisher={Spe}
}

@article{cao2022bayesian,
  title={Bayesian calibration for large-scale fluid structure interaction problems under embedded/immersed boundary framework},
  author={Cao, Shunxiang and Zhengyu Huang, Daniel},
  journal={International Journal for Numerical Methods in Engineering},
  volume={123},
  number={8},
  pages={1791--1812},
  year={2022},
  publisher={Wiley Online Library}
}

@article{dashti2013map,
  title={{MAP} estimators and their consistency in {B}ayesian nonparametric inverse problems},
  author={Dashti, Masoumeh and Law, Kody JH and Stuart, Andrew M and Voss, Jochen},
  journal={Inverse Problems},
  volume={29},
  number={9},
  pages={095017},
  year={2013},
  publisher={IOP Publishing}
}

@article{armbruster2022stabilization,
  title={A stabilization of a continuous limit of the ensemble Kalman inversion},
  author={Armbruster, Dieter and Herty, Michael and Visconti, Giuseppe},
  journal={SIAM Journal on Numerical Analysis},
  volume={60},
  number={3},
  pages={1494--1515},
  year={2022},
  publisher={SIAM}
}

@article{al2024non,
  title={Non-asymptotic analysis of ensemble Kalman updates: effective dimension and localization},
  author={Al-Ghattas, Omar and Sanz-Alonso, Daniel},
  journal={Information and Inference: A Journal of the IMA},
  volume={13},
  number={1},
  pages={iaad043},
  year={2024},
  publisher={Oxford University Press}
}

@article{tong2023localized,
  title={Localized ensemble Kalman inversion},
  author={Tong, Xin T and Morzfeld, Matthias},
  journal={Inverse Problems},
  volume={39},
  number={6},
  pages={064002},
  year={2023},
  publisher={IOP Publishing}
}

@article{hendricks2008real,
  title={Real-time groundwater flow modeling with the ensemble Kalman filter: Joint estimation of states and parameters and the filter inbreeding problem},
  author={Hendricks Franssen, Harrie-Jan and Kinzelbach, W},
  journal={Water resources research},
  volume={44},
  number={9},
  year={2008},
  publisher={Wiley Online Library}
}

@article{aanonsen2009ensemble,
  title={The ensemble Kalman filter in reservoir engineering—a review},
  author={Aanonsen, Sigurd I and N{\oe}vdal, Geir and Oliver, Dean S and Reynolds, Albert C and Vall{\`e}s, Brice},
  journal={Spe Journal},
  volume={14},
  number={03},
  pages={393--412},
  year={2009},
  publisher={SPE}
}

@article{qian2025fundamental, title={The fundamental subspaces of ensemble Kalman inversion}, author={Qian, Elizabeth and Beattie, Christopher}, journal={SIAM Review}, volume={67}, number={4}, pages={771--798}, year={2025}, publisher={SIAM} }

@article{blomker2019well,
  title={Well posedness and convergence analysis of the ensemble Kalman inversion},
  author={Bl{\"o}mker, Dirk and Schillings, Claudia and Wacker, Philipp and Weissmann, Simon},
  journal={Inverse Problems},
  volume={35},
  number={8},
  pages={085007},
  year={2019},
  publisher={IOP Publishing}
}

@article{iglesias2021adaptive,
  title={Adaptive regularisation for ensemble {K}alman inversion},
  author={Iglesias, Marco and Yang, Yuchen},
  journal={Inverse Problems},
  volume={37},
  number={2},
  pages={025008},
  year={2021},
  publisher={IOP Publishing}
}

@article{che2025stable,
  title={Stable derivative free Gaussian mixture variational inference for Bayesian inverse problems},
  author={Che, Baojun and Chen, Yifan and Huan, Zhenghao and Huang, Daniel Zhengyu and Wang, Weijie},
  journal={SIAM Journal on Scientific Computing},
  volume={47},
  number={5},
  pages={A2583--A2608},
  year={2025},
  publisher={SIAM}
}

@article{che2026adaptive, title={Adaptive Exponential Integration for Stable Gaussian Mixture Black-Box Variational Inference}, author={Che, Baojun and Chen, Yifan and Huang, Daniel Zhengyu and Mao, Xinying and Wang, Weijie}, journal={arXiv preprint arXiv:2601.14855}, year={2026} }

@article{yan2019adaptive,
  title={An adaptive multifidelity PC-based ensemble Kalman inversion for inverse problems},
  author={Yan, Liang and Zhou, Tao},
  journal={International Journal for Uncertainty Quantification},
  volume={9},
  number={3},
  year={2019},
  publisher={Begel House Inc.}
}

@article{hansen2003reducing,
  title={Reducing the time complexity of the derandomized evolution strategy with covariance matrix adaptation (CMA-ES)},
  author={Hansen, Nikolaus and M{\"u}ller, Sibylle D and Koumoutsakos, Petros},
  journal={Evolutionary computation},
  volume={11},
  number={1},
  pages={1--18},
  year={2003},
  publisher={MIT Press}
}

@article{yang2026polarized,
  title={Polarized Ensemble Kalman Inversion: A Derivative-Free Method for Multimodal Inverse Problems with Non-asymptotic Analysis},
  author={Yang, Huchen and Wu, Jin-Long},
  journal={arXiv preprint arXiv:2610.31907},
  year={2026}
}

@article{gao2024adaptive,
  title={Adaptive operator learning for infinite-dimensional Bayesian inverse problems},
  author={Gao, Zhiwei and Yan, Liang and Zhou, Tao},
  journal={SIAM/ASA Journal on Uncertainty Quantification},
  volume={12},
  number={4},
  pages={1389--1423},
  year={2024},
  publisher={SIAM}
}

@article{chen2024efficient,
  title={Efficient, multimodal, and derivative-free Bayesian inference with Fisher--Rao gradient flows},
  author={Chen, Yifan and Zhengyu Huang, Daniel and Huang, Jiaoyang and Reich, Sebastian and Stuart, Andrew M},
  journal={Inverse Problems},
  volume={40},
  number={12},
  pages={125001},
  year={2024},
  publisher={IOP Publishing}
}

@article { held1985pseudomomentum,
      author = "Isaac M.  Held",
      title = "Pseudomomentum and the Orthogonality of Modes in Shear Flows",
      journal = "Journal of Atmospheric Sciences",
      year = "1985",
      publisher = "American Meteorological Society",
      address = "Boston MA, USA",
      volume = "42",
      number = "21",
      doi = "10.1175/1520-0469(1985)042<2280:PATOOM>2.0.CO;2",
      pages=      "2280 - 2288",
      url = "https://journals.ametsoc.org/view/journals/atsc/42/21/1520-0469_1985_042_2280_patoom_2_0_co_2.xml"
}

@article {held1987linear,
      author = "Isaac M.  Held and Peter J.  Phillips",
      title = "Linear and Nonlinear Barotropic Decay on the Sphere",
      journal = "Journal of Atmospheric Sciences",
      year = "1987",
      publisher = "American Meteorological Society",
      address = "Boston MA, USA",
      volume = "44",
      number = "1",
      doi = "10.1175/1520-0469(1987)044<0200:LANBDO>2.0.CO;2",
      pages=      "200 - 207",
      url = "https://journals.ametsoc.org/view/journals/atsc/44/1/1520-0469_1987_044_0200_lanbdo_2_0_co_2.xml"
}
\end{document}